\documentclass[12pt, reqno]{amsart}

\author[M.~Caprio]{Michele Caprio and Sayan Mukherjee}
\address{PRECISE Center, Dept. of Computer and Information  Science, 
University of Pennsylvania, 3330 Walnut Street, Philadelphia, PA 19104}
\email{caprio@seas.upenn.edu}
\urladdr{\url{https://mc6034.wixsite.com/caprio}}  
\address{Center for Scalable Data Analytics and Artificial Intelligence, Universität Leipzig, Humboldtstraße 25, Leipzig, Germany 04105 and the Max Planck Institute for Mathematics in the Sciences, Inselstraße 22
04103 Leipzig
Germany; Dept. of Statistical Science, Mathematics, Computer Science, and Biostatistics \& Bioinformatics, Duke University, Durham, NC 27708, USA}
\email{sayan.mukherjee@mis.mpg.de}
\urladdr{\url{https://sayanmuk.github.io/}}


\title{Ergodic Theorems for Dynamic Imprecise Probability Kinematics}

\usepackage[T1]{fontenc}
\usepackage{amsmath}
\usepackage{amssymb}
\usepackage{amsthm}
\usepackage[left=2.5cm, right=2.5cm, top=3cm]{geometry}
\usepackage{hyperref}
\usepackage{algorithm}
\usepackage{algpseudocode}
\usepackage{bbm}
\usepackage{tikz}
\usepackage{caption}
\usepackage{subcaption}
\usepackage{fancyhdr}
\usepackage[inline]{enumitem}
\usepackage{comment}
\usepackage{nicefrac}
\usepackage{bm}
\usepackage{mathrsfs}
\usepackage{graphicx}
\usepackage[utf8]{inputenc}
\usepackage{cancel}
\usepackage{mathtools}

\newcommand{\vertiii}[1]{{\left\vert\kern-0.25ex\left\vert\kern-0.25ex\left\vert #1 
    \right\vert\kern-0.25ex\right\vert\kern-0.25ex\right\vert}}
		
\AtBeginDocument{%
   \def\MR#1{}
}

\newtheorem{thm}{Theorem}[section]
\theoremstyle{definition} 
\newtheorem{defi}[thm]{Definition}%[section]
\let\olddefi\defi
\renewcommand{\defi}{\olddefi\normalfont}
\newtheorem{rmk}[thm]{Remark}
\let\oldrmk\rmk
\renewcommand{\rmk}{\oldrmk\normalfont}

\newtheorem{theorem}{Theorem}%[section]
\newtheorem{lemma}[theorem]{Lemma}

\newtheorem{corollary}[theorem]{Corollary}

\newtheorem{claim}[theorem]{Claim}

\hypersetup{
    pdftitle={prove},
    pdfauthor={autori},
    pdfmenubar=false,
    pdffitwindow=true,
    pdfstartview=FitH,
    colorlinks=true,
    linkcolor=blue,
    citecolor=green,
    urlcolor=cyan
}

\providecommand{\MR}[1]{}

\providecommand{\MR}{\relax\ifhmode\unskip\space\fi MR }

\providecommand{\href}[2]{#2}

\begin{document}

\maketitle
\thispagestyle{empty}

\begin{abstract}
We formulate an ergodic theory for the (almost sure) limit $\mathcal{P}^\text{co}_{\tilde{\mathcal{E}}}$ of a sequence $(\mathcal{P}^\text{co}_{\mathcal{E}_n})$ of successive dynamic imprecise probability kinematics (DIPK, introduced in \cite{prob.kin}) updates of a set $\mathcal{P}^\text{co}_{\mathcal{E}_0}$ representing the initial beliefs of an agent. As a consequence, we formulate a strong law of large numbers.
\end{abstract}

\section{Introduction}\label{intro}
In \cite{prob.kin} the authors introduce a procedure called dynamic imprecise probability kinematics (DIPK) to update an agent's opinions on the elements of sigma-algebra $\mathcal{F}=2^\Omega$ of a state space of interest $\Omega$ (at most countable) in the presence of ambiguity and partial information. The former refers to a situation in which a single probability measure is not enough to encapsulate the agent's initial beliefs. To account for this, the agent specifies a (closed and convex) set of probability measures $\mathcal{P}^\text{co}_{\mathcal{E}_0}$ representing their initial state of knowledge. 
%a finite number of plausible probabilities $P_1,\ldots,P_k$ representing their initial state of knowledge, which induce a \textit{(finitely generated) credal set}, that for DIPK is their convex hull $\mathcal{P}^\text{co}_{\mathcal{E}_0}=\text{Conv}(P_1,\ldots,P_k)$. 
%a (closed and convex) set of probability measures. representing the agent's initial for $\mathcal{P}^{\text{co}}_{\tilde{\mathcal{E}}}$ beliefs the (time)
Partial information means that the agent cannot collect crisp evidence; rather, they gather  information whose nature is probabilistic.

In this work, we use tools from  dynamical systems theory to develop an ergodic theory for the limiting set of probabilities $\mathcal{P}^\text{co}_{\tilde{\mathcal{E}}}$ of the sequence $(\mathcal{P}^\text{co}_{\mathcal{E}_n})$ of DIPK updates of set $\mathcal{P}^\text{co}_{\mathcal{E}_0}$. The ergodic theorems that we present in this paper are instrumental to increase the applicability of DIPK; for example, they underpin generalizations of classical MCMC procedures that allow for DIPK updating. Such methods will be the subject of future studies. 

Our main result is Theorem \ref{ergo1}. We show that, given a (bounded) functional $f$ on $\Omega$ and an operator $T:\Omega \rightarrow\Omega$, the average $\frac{1}{k} \sum_{j=1}^k f(T^{j-1}(\omega))$ may not converge, but -- under some regularity assumptions -- will almost surely be eventually contained in the interval 
$$\left[\inf_{P_{\tilde{\mathcal{E}}}\in\mathcal{P}^{\text{co}}_{\tilde{\mathcal{E}}}}\sum_{\omega \in \Omega} f^\star(\omega) {P}_{\tilde{\mathcal{E}}}(\{\omega\}) , \sup_{P_{\tilde{\mathcal{E}}}\in\mathcal{P}^{\text{co}}_{\tilde{\mathcal{E}}}}\sum_{\omega \in \Omega} f^\star(\omega) {P}_{\tilde{\mathcal{E}}}(\{\omega\})\right],$$ 
for a well defined function $f^\star$ on $\Omega$. This means that the limit infimum and the limit supremum of $\frac{1}{k} \sum_{j=1}^k f(T^{j-1}(\omega))$ coincide and belong almost surely to the aforementioned interval. Their endpoints are given by the infimum and the supremum of the expectation of $f^\star$ with respect to the probability measures in  $\mathcal{P}^\text{co}_{\tilde{\mathcal{E}}}$.  
%lowest and the highest expected values of $f^\star$ with respect to the elements of $\mathcal{P}^\text{co}_{\tilde{\mathcal{E}}}$. 
%(a version of) the expected values of $f^\star$ with respect to the lower and upper probability measures associated with $\mathcal{P}^\text{co}_{\tilde{\mathcal{E}}}$. 
%, the limit set of the sequence $(\mathcal{P}^\text{co}_{{\mathcal{E}}_n})$ of DIPK updates of set $\mathcal{P}^\text{co}_{\mathcal{E}_0}$ representing the initial beliefs of the agent. 
Notice that here almost surely means that the lower probability $\underline{P}_{\tilde{\mathcal{E}}}$ that the event happens is $1$.
\subsection{Structure of the paper}
In an effort to make the paper self-contained, we give the needed preliminary notions in section \ref{prelim}. Section \ref{ergo.result} contains our results and section \ref{concl} concludes our work. Appendix \ref{app_lemmas} consists of technical lemmas, appendix \ref{dpk} deals with the dynamic (precise) probability kinematics (DPK) updating procedure, while proofs can be found in appendix \ref{proofs}.

\section{Preliminaries}\label{prelim}
In this section, we give the preliminaries needed to understand the main results of the paper. 
\subsection{Non-additive notions}\label{NAD_imp}
We begin with the non-additive notions needed for our analysis; we follow \cite{cerreia}. 
%, we introduce the following.  
Given a measurable space $(\Omega,\mathcal{F})$, where $\Omega \neq \emptyset$ and $\mathcal{F}$ is a sigma-algebra of subsets of $\Omega$, we say that a set function $\nu:\mathcal{F} \rightarrow [0,1]$ is a Choquet capacity if $\nu(\emptyset)=0$, $\nu(\Omega)=1$, and $\nu(A) \leq \nu(B)$ for all $A,B \in \mathcal{F}$ such that $A \subset B$. Call $\Delta(\Omega,\mathcal{F})$ the set of probability measures on $(\Omega,\mathcal{F})$. Then, we say that a Choquet capacity $\nu:\mathcal{F} \rightarrow [0,1]$ is
\begin{itemize}
\item[(i)] \textit{convex} if $\nu(A \cup B) + \nu(A \cap B) \geq \nu(A) +\nu(B)$, for all $A,B \in \mathcal{F}$;
\item[(ii)] \textit{additive} if $\nu(A \cup B)=\nu(A) +\nu(B)$, for all disjoint $A,B \in \mathcal{F}$;
\item[(iii)] \textit{continuous} if $\lim_{n \rightarrow \infty} \nu(A_n)=\nu(A)$ whenever either $A_n \uparrow A$ or $A_n \downarrow A$;
\item[(iv)] \textit{continuous at $\Omega$} if $\lim_{n \rightarrow \infty} \nu(A_n)=\nu(\Omega)$ whenever $A_n \uparrow \Omega$;
\item[(v)] a \textit{probability measure} if it is an additive Choquet capacity which is continuous at $\Omega$;
\item[(vi)] a \textit{lower probability measure} if there exists a set $\mathcal{P} \subset \Delta(\Omega,\mathcal{F})$ such that 
$$\nu(A)=\inf_{P \in \mathcal{P}}P(A), \quad \text{for all } A \in \mathcal{F}.$$
\end{itemize}
%\end{definition}
The \textit{upper probability measure} $\overline{\nu}$ associated with $\mathcal{P}$ is defined as the conjugate to $\nu$, that is, $\overline{\nu}(A):=1-{\nu}(A^c)=\sup_{P \in \mathcal{P}} P(A)$, for all $A \in \mathcal{F}$. A generic lower probability $\nu$ completely characterizes the set of probability measures that setwise dominate $\nu$, called the core of $\nu$
\begin{align*}
   \text{core}(\nu)&:=\{P \in \Delta(\Omega,\mathcal{F}): P(A) \geq \nu(A), \text{ for all } A \in \mathcal{F}\}.
\end{align*}
By completely characterize, we mean that it is enough to know $\nu$ to retrieve all the elements in the core. The core is convex \cite[Section 2.2]{marinacci_ambig} and weak$^\star$-compact \cite[Proposition 3]{marinacci_ambig}.\footnote{In the weak$^\star$ topology, a net $(P_\alpha)_{\alpha \in I}$ converges to $P$ if and only if $P_\alpha(A) \rightarrow P(A)$, for all $A \in \mathcal{F}$.}

We say that a generic lower probability $\nu$ is
\begin{itemize}
    \item[($\alpha$)] \textit{($T$-)invariant} if, for all $A \in \mathcal{F}$,
$$\nu(A)=\nu(T^{-1}(A)) \text{.}$$
We then call $\mathcal{I} \subset \Delta(\Omega,\mathcal{F})$ the set of ($T$-){invariant probability measures}, that is, 
$$\mathcal{I}:=\left\lbrace{P \in \Delta(\Omega,\mathcal{F}) : P(A)=P(T^{-1}(A)) \text{, for all } A \in \mathcal{F}}\right\rbrace.$$
We call $\mathcal{G} \in \mathcal{F}$ the set of all ($T$-){invariant events} of $\mathcal{F}$, that is, 
$$\mathcal{G}:=\left\lbrace{A \in \mathcal{F} : T^{-1}(A)=A}\right\rbrace.$$
    \item[($\beta$)] \textit{ergodic} if and only if $\nu(\mathcal{G})=\{0,1\}$, that is, $\nu$ assigns value $0$ or $1$ to all the elements of $\mathcal{G}$.
    \item[($\gamma$)] \textit{strongly invariant} if and only if for every $A \in \mathcal{F}$, 
$$\nu(A\setminus T^{-1}(A))=\overline{\nu}(T^{-1}(A) \setminus A) \quad \text{and} \quad {\nu}(T^{-1}(A) \setminus A)=\overline{\nu}(A\setminus T^{-1}(A));$$ 
    \item[($\delta$)] \textit{functionally invariant} if and only if $\mathcal{M} \subset \mathcal{I}$, where $\mathcal{M} \subset \Delta(\Omega,\mathcal{F})$ is the set for which $\nu(A)=\inf_{P \in \mathcal{M}} P(A)$, for all $A \in \mathcal{F}$.
 \end{itemize}   
Finally, given a generic lower probability $\nu$ and a generic function $f\in B(\Omega,\mathcal{F})$, the set of bounded and $\mathcal{F}$-measurable functionals on $\Omega$, we define the \textit{Choquet integral} as follows
$$\int_\Omega f \text{ d}\nu:=\int_0^\infty \nu\left(\{\omega\in\Omega : f(\omega)\geq t\} \right) \text{ d}t + \int_{-\infty}^0 \left[\nu\left(\{\omega\in\Omega : f(\omega)\geq t\} \right) -\nu(\Omega) \right] \text{ d}t,$$
where the right hand side integrals are (improper) Riemann integrals. If $\nu$ is additive, then the Choquet integral reduces to the standard additive integral. 

Consider now a set of probability measures $\mathcal{P}$ on a generic measurable space $(\Omega,\mathcal{F})$ and its lower envelope $\nu$ -- that is, $\nu(A)=\inf_{P\in\mathcal{P}}P(A)$, for all $A \in\mathcal{F}$. Then, \cite{walley_report} points out that given the Choquet integral of any function $f\in B(\Omega,\mathcal{F})$ with respect to $\nu$ is a lower bound to the infimum of the expectation of $f$ with respect to the probability measures in  $\mathcal{P}$. In formulas, 
\begin{equation}\label{cvx_bd}
    \int_\Omega f \text{ d} \nu \leq \inf_{P\in\mathcal{P}} \int_\Omega f \text{ d} P = \inf_{P\in\mathcal{P}} \mathbb{E}_P(f), \quad \text{for all }  f\in B(\Omega,\mathcal{F}).
\end{equation}
In \cite{walley_report} the author also shows that the inequality in \eqref{cvx_bd} holds with equality if and only if $\nu$ is convex.

\subsection{Ergodic theory}\label{why_ergo}
Ergodic theory is a branch of mathematics that studies the long-term average behavior
of complex dynamical systems; it was  first introduced 
%by the Austrian physicist Ludwig Boltzmann 
in \cite{boltzmann}. Working with gases, the author suggested that the spatial average values giving rise to macroscopic features also arose as averages over time of observable quantities that could be calculated from microscopic states. Hence, what can be considered the ``ergodic mantra'': \textit{space average equals time average}. The best-known ergodic theorem is arguably Birkhoff's ergodic theorem (cf. \cite{cornfeld}).
\begin{theorem}\label{birk}\textbf{(Birkhoff)}
If we have a probability space $(\Omega,\mathcal{F},P)$, a measurable self-map $T:\Omega \rightarrow \Omega$ such that $P(A)=P(T^{-1}(A))$, for all $A \in \mathcal{F}$, and a measurable functional $f$ on $\Omega$, then the limit as $n \rightarrow \infty$ of the time average $\frac{1}{n} \sum_{j=1}^n f(T^{j-1}(\omega))$ exists and it is equal to the space average $\frac{1}{P(\Omega)} \int_\Omega f \text{d}P=\int_\Omega f \text{d}P$, $P -$a.s.
\end{theorem}

%is especially meaningful because it
This result characterizes the behavior of the orbit of operator $T$ over a long time period. In particular, the time average $$\frac{1}{n} \sum_{j=1}^n f(T^{j-1}(\omega))$$ of $f$ will almost surely converge to $\mathbb{E}_P(f)=\int_\Omega f \text{ d}P$ as the time horizon $n$ approaches infinity, where ``almost surely'' means that the probability that it does not happen is zero.

The importance of ergodic theory for computer scientists and statisticians is well documented in many works. For example, in \cite{berry} the authors combine ideas from the theory of dynamical systems with learning theory, providing an effective route to data-driven models of complex systems. They obtain refinable predictions as the amount of training data increases, and physical interpretability through discovery of coherent patterns around which the dynamics is organized. In \cite{zhang} the author proposes a generalization of ergodic measure preserving flow (EMPF), an optimisation-based
inference method using ergodic results that overcomes the biasedness limitations of both Markov chain Monte Carlo (MCMC) and variational inference (VI). Such generalization, called  ergodic inference, is necessary because of the lack of theoretical proof
of the validity of EMPF. In \cite{arbabi} the authors establish the convergence of a class of numerical algorithms, known as dynamic mode decomposition (DMD), for computation of the eigenvalues and eigenfunctions of the infinite dimensional Koopman operator. Koopman operator theory is an alternative formulation of dynamical systems theory which
provides a versatile framework for data-driven study of high-dimensional nonlinear systems. Their work rely on the assumption that the underlying dynamical system is ergodic. In \cite{vose} the author points out that genetic algorithms are strongly related to dynamical systems. Ergodicity of such systems corresponds to an important property, called asymptotic correctness, roughly guaranteeing to eventually explore the whole solution space.

Ergodic theorems for imprecise probabilities (IP) have been studied in the context of imprecise Markov chains in \cite{de_bock, de_cooman}, capacity-preserving $\mathbb{Z}^d_+$-actions \cite{wu}, and lower probabilities in \cite{cerreia}. In this latter the authors work with an uncountable state space $\Omega$. In section \ref{ergo.result} we provide similar results in the context of dynamic imprecise probability kinematics. 

Three other prominent papers in the IP literature concerning ergodicity are \cite{hermans,de_cooman3,tjones}. In the first one, the authors study the ergodicity of upper transition operators (UTOs), that are bounded, subadditive, and non-negatively homogeneous transformations of finite-dimensional linear spaces. They were introduced in \cite[section 3]{de_cooman} when describing imprecise Markov chains, random processes where prior and transition beliefs are described in terms of coherent upper previsions (CHPs) \cite[Section 2.3.3]{walley}. In \cite[Section 2]{hermans}, the authors relate UTOs and CHPs: they show how upper transition operators can be seen as the Cartesian product of the coherent upper conditional previsions over all states. The main  differences between \cite{hermans} and the present paper are the following. We do not consider UTOs; rather, in section \ref{ergo.result} we work with a generic $\mathcal{F}\backslash\mathcal{F}$-measurable transformation on $\Omega$, that is, a generic operator $T$. In addition, while in \cite{hermans} UTO $T$ operates on the finite-dimensional linear space $\mathcal{L}(\Omega)$ of real-valued maps on a finite nonempty state space $\Omega$, in this paper $T$ is defined on $\Omega$ (finite or countable). This means that while in our work results concern (bounded measurable) functional $f$ evaluated at (the $j$-th orbit of) operator $T$ -- written $f(T^j(\omega))$, for some $\omega\in\Omega$ --, in \cite{hermans} results concern (the $j$-th orbit of) UTO $T$ evaluated at functional $f$ -- written $T^j(f(\omega))$, for some $\omega\in\Omega$.\footnote{We believe our approach is more consistent with classical ergodic theory, see e.g. the formulation of Birkhoff's ergodic theorem in Theorem \ref{birk}. The authors of \cite{hermans} seem to agree with us, see \cite[Section 3]{hermans}.} It is also worth mentioning that the object of study of \cite{hermans} is the ergodicity of a UTO, while in our work it is the ergodicity of lower probability $\underline{P}_{\tilde{\mathcal{E}}}$.

In \cite{de_cooman3}, the authors study joint lower and upper expectations in imprecise probability trees (IPTs), in terms of the sub- and supermartingales associated with such trees; in particular, they derive a law of iterated expectations. IPTs can be seen as discrete-time stochastic processes with finite state spaces and imprecise transition probabilities (they belong to a credal set, that is, a set of probability measures). They then focus on imprecise Markov chains, and study many of their properties. Notably, and more significantly in the context of this paper, they prove a game-theoretic version of the strong law of large numbers for submartingale differences in IPTs, and use it to derive point-wise ergodic theorems for imprecise Markov chains, involving (bounded) functions of a finite number of states. Paraphrasing what the authors state in \cite[Section 1]{de_cooman3}, the ergodic results provided in our work are at once more general and more restricted than those in \cite{de_cooman3}. In particular, in this paper the context is not restricted to shift invariance in Markov chains, but extends to invariance under arbitrary transformations on at most countable sample spaces. On the other hand, the assumptions required in section \ref{ergo.result} for our ergodic results are somehow stronger than those required in \cite[Section 9]{de_cooman3} (see also Lemmas \ref{suff_cond_invariance}, \ref{suff_cond_inv_ergo}, \ref{suff_cond_stricter_than_before}, and \ref{suff_cond_stricter_than_before2}, in appendix \ref{app_lemmas}). Our results, then, cannot be obtained as a special case of those in \cite{de_cooman3} and vice versa.

Finally, in \cite{tjones} the authors study the limit behavior of upper and lower bounds on expected time averages (ETAs) in imprecise Markov chains. They find a necessary and sufficient condition under which these upper and lower bounds converge as time goes to infinity to limit values that do not depend on the initial state. Given that this condition is considerably weaker than that needed for ergodic behavior, they refer to their result as weak ergodicity (WE). The characterization of WE, as well as the values of  upper and lower ETAs do not depend on the types of independence assumed (introduced in \cite[Section 1]{tjones}). 
%; this type of robustness is not is not exhibited by the notion of ergodicity and the related inferences of limit upper and lower expectations. 
They conclude by showing that directly using upper and lower expected time averages improves the information about the limit behavior of time averages with respect to the use of limit upper and lower expectations. \cite{tjones} is a very profound and general article; its results, though, do not subsume ours because, as pointed out when comparing the present paper with \cite{hermans}, we consider a different context, namely, we do not work with imprecise Markov chains.
%those in the present work because of two main differences, already pointed our when comparing \cite{hermans} and this paper: we do not work with imprecise Markov chains and we do not require the state space $\Omega$ to be finite. 
%The authors derive numerous properties for their joint lower and upper expectations, and in particular a version of the tower rule. 
\subsection{Dynamic imprecise probability kinematics}
 In \cite{prob.kin}, the authors require the state space $\Omega$ to be at most countable, and the sigma-algebra on $\Omega$ to be the power set $\mathcal{F}=2^\Omega$. The latter assumption is made to work with the richest possible sigma-algebra; all the results in \cite{prob.kin} and in the present paper still hold if $\mathcal{F}$ is not the power set. $\Omega$ is assumed at most countable for simplicity (the authors wanted to focus on the updating mechanism and not on measure-theoretic complications).
 
 % and in view of future developments, as explained in a few lines
 
 DIPK updating is described as follows. First, the authors prescribe the agent to specify a finite set of plausible probabilities $\mathcal{P}=\{P_1,\ldots,P_k\}$ capturing the ambiguity they face, and to compute the lower probability associated with it. The core of such lower probability, denoted by $\mathcal{P}^\text{co}_{\mathcal{E}_0}$, represents the agent's initial beliefs. The authors work with the core because it is completely characterized by the lower probability and because in general the core is a superset of the convex hull of $P_1,\ldots,P_k$ \cite[Remark 14]{prob.kin}, thus capturing a higher level of uncertainty. The authors (tacitly) assume that the extrema of $\mathcal{P}^\text{co}_{\mathcal{E}_0}$, that is, the elements of the core that cannot be written as a convex combination of other elements, are finite; they do so for the following reasons. 
 \begin{itemize}
     \item If $\Omega$ is finite, we can see probability measures as vectors in the unit simplex $\Delta^{d-1}$ of $\mathbb{R}^d$, where $d=\#\Omega<\infty$ and $\#$ represents the cardinality operator. So the core of a lower probability will be a (closed) convex subset of $\Delta^{d-1}$, which can be approximated arbitrarily well by a polytope having finitely many vertices \cite{bronstein}.\footnote{Here ``approximated arbitrarily well'' means that some distance between the convex set and the polytope, e.g. the Hausdorff metric, can be made arbitrarily small.} The polytope with finitely many vertices is the geometric representation of a closed and convex set of probabilities having finitely many extrema.
     \item If $\Omega$ is countable, the assumption is stronger, and has mainly a computational motivation. It corresponds to the agent specifying a \textit{finitely generated credal set} -- that is, the convex hull of finitely many probability measures -- that is (possibly) a superset of the convex hull of $P_1,\ldots,P_k$, $\text{Conv}(P_1,\ldots,P_k)$. 
 \end{itemize}
 In the future, the authors plan to generalize DIPK by working with lower previsions of random variables instead of just focusing on lower probabilities.\footnote{A (bounded) random variable $X$ can be seen as a gamble yielding uncertain rewards, expressed in units of utility \cite{walley}. The lower prevision $\underline{P}(X)$ of $X$ is the highest price that the agent accepts to buy gamble $X$, that is, the highest $\alpha$ so that the agent accepts uncertain reward $X-\alpha$. If we have a set $\mathcal{P}$ of probability measures, then $\underline{P}(X)=\inf_{P\in\mathcal{P}} \mathbb{E}_P(X)$.} As Lemma \ref{core_conv} in appendix \ref{app_lemmas} shows, in that case the core (appropriately redefined) and $\text{Conv}(P_1,\ldots,P_k)$ coincide, so the assumption that they are equal is automatically verified.
 
 To update their beliefs, the agent computes the dynamic (precise) probability kinematics (DPK) update  of the extrema of the core. DPK updating is studied in detail in \cite[Sections 3-5]{prob.kin}; we introduce it briefly in appendix \ref{dpk}. Their updated beliefs are represented by the convex hull of the updated extrema, which coincides with the core of the updated lower probability -- that is, it coincides with updating every element in the core and then taking the lower envelope -- by the following theorem, reported in \cite{walley}.
\begin{theorem}\label{choq_rel}
Suppose $\text{core}(\nu) \neq \emptyset$. Then, the following holds.
\begin{itemize}
\item[(a)] The set of the extrema of $\text{core}(\nu)$ is nonempty, in symbols $\text{ex}(\text{core}(\nu)) \neq \emptyset$.
\item[(b)] $\text{core}(\nu)$ is the closure in the weak$^\star$ topology of the convex hull of $\text{ex}(\text{core}(\nu))$.
\item[(c)] If $\nu(A)=\inf_{P \in \text{core}(\nu)} P(A)$, for all $A \in \mathcal{F}$, then $\nu(A)=\inf_{P \in \text{ex}(\text{core}(\nu))} P(A)$, for all $A \in \mathcal{F}$.
\end{itemize}
\end{theorem}

Repeating this procedure gives us the sequence $(\mathcal{P}^{\text{co}}_{\mathcal{E}_t})$ of successive DIPK updates of $\mathcal{P}^\text{co}_{\mathcal{E}_0}$. The elements of the sequence are indexed by $\mathcal{E}_t$, a partition of $\Omega$ induced by the available data at time $t$ \cite[Sections 3-5]{prob.kin}. Under mild regularity conditions, $(\mathcal{P}^{\text{co}}_{\mathcal{E}_t})$ converges to a set denoted as  $\mathcal{P}^\text{co}_{\tilde{\mathcal{E}}}$ \cite[Proposition 12]{prob.kin}.\footnote{The regularity conditions are the following: the data points that drive the updating procedure are sampled independently, random variable $X$ introduced in appendix \ref{dpk} has a finite first moment, and coefficient $\beta(n)$ introduced in appendix \ref{dpk} is $o(1/n)$.} In the remainder of the paper, we assume these regularity conditions are satisfied. The results presented in this work are based on lower probability 
%$\underline{P}_{\tilde{\mathcal{E}}}$, 
%defined on measurable space $(\Omega,\mathcal{F})$, with $\Omega$ at most countable and $\mathcal{F}=2^\Omega$,
$\underline{P}_{\tilde{\mathcal{E}}}(A)=\inf_{P_{\tilde{\mathcal{E}}}\in \mathcal{P}^\text{co}_{\tilde{\mathcal{E}}}} P_{\tilde{\mathcal{E}}}(A)$, for all $A\in\mathcal{F}$, because $\underline{P}_{\tilde{\mathcal{E}}}$ completely characterizes $\mathcal{P}^\text{co}_{\tilde{\mathcal{E}}}$.
%, as shown in \cite[Section 6]{prob.kin}.

%Recall the following definition, given in 
%Following \cite{cerreia}, we introduce some concepts that are crucial for understanding the results in the paper.
%\begin{definition}\label{choq_cap_gen}
 
%That is, it is enough to know $\underline{P}_{\tilde{\mathcal{E}}}$ to be able to specify all the elements in the set.

%\subsection{Why ergodic theory?}\label{why_ergo}

\section{Ergodic theory for $\mathcal{P}^\text{co}_{\tilde{\mathcal{E}}}$}\label{ergo.result}

Recall that in the DIPK procedure described in \cite{prob.kin}, $\Omega$ is assumed finite or countable. Let $T: \Omega \rightarrow \Omega$ be an $\mathcal{F} \backslash \mathcal{F}$-measurable transformation (this corresponds to the ergodic operator in classical ergodic theory) that explores all the state space, that is, for all $\omega \in \Omega$,
$$\bigcup_{j \in \mathbb{N}} T^{j-1}(\omega)=\Omega.$$
Let also $f:\Omega \rightarrow \mathbb{R}$ belong to $B(\Omega,\mathcal{F})$, the set of bounded and $\mathcal{F}$-measurable functionals on $\Omega$. Then, under some regularity conditions, the limit of the empirical average $\frac{1}{k} \sum_{j=1}^{k} f(T^{j-1} (\omega))$ as $k \rightarrow \infty$ exists $\underline{P}_{\tilde{\mathcal{E}}} -$a.s., and it belongs to the interval generated by the infimum and the supremum of the space average of $f^\star$ with respect to the probability measures in  $\mathcal{P}^\text{co}_{\tilde{\mathcal{E}}}$,
%(a version of) the space averages taken with respect to the boundary elements of $\mathcal{P}^\text{co}_{\tilde{\mathcal{E}}}$, lowest and highest space averages of $f^\star$ taken with respect to 
$$\mathcal{A}(\omega):=\left[\inf_{P_{\tilde{\mathcal{E}}}\in\mathcal{P}^{\text{co}}_{\tilde{\mathcal{E}}}}\sum_{\omega \in \Omega} f^\star(\omega) {P}_{\tilde{\mathcal{E}}}(\{\omega\}) , \sup_{P_{\tilde{\mathcal{E}}}\in\mathcal{P}^{\text{co}}_{\tilde{\mathcal{E}}}}\sum_{\omega \in \Omega} f^\star(\omega) {P}_{\tilde{\mathcal{E}}}(\{\omega\})\right],$$
$\underline{P}_{\tilde{\mathcal{E}}} -$a.s., for a well defined functional $f^\star$. That is, 
$$\underline{P}_{\tilde{\mathcal{E}}} \left(\left\lbrace{\omega \in \Omega : \lim_{k\rightarrow\infty} \frac{1}{k} \sum_{j=1}^{k} f(T^{j-1} (\omega)) \in \mathcal{A}(\omega)}\right\rbrace \right)=1.$$

%We now present three lemmas that 
The following is the main contribution of the paper. It can be seen as the DIPK version of \cite[Theorem 2]{cerreia} with an extra result, namely equation \eqref{ergo.result_sum}.
\begin{theorem}\label{ergo1}
If $\underline{P}_{\tilde{\mathcal{E}}}$ is invariant, then for all $f \in B(\Omega,\mathcal{F})$, there exists $f^\star \in B(\Omega,\mathcal{G})$ -- that is, there exists a bounded and $\mathcal{G}$-measurable functional $f^\star$ on $\Omega$ -- such that
\begin{align}\label{lim.emp.avg}
\lim\limits_{k \rightarrow \infty} \frac{1}{k} \sum_{j=1}^{k} f(T^{j-1} (\omega)) = {f}^\star (\omega) \quad \underline{P}_{\tilde{\mathcal{E}}}-a.s.
\end{align}
If in addition $\underline{P}_{\tilde{\mathcal{E}}}$ is ergodic, then
\begin{equation}\label{ergo.result_int}
\int_\Omega f^\star \text{ d}\underline{P}_{\tilde{\mathcal{E}}} \leq \lim\limits_{k \rightarrow \infty} \frac{1}{k} \sum_{j=1}^{k} f(T^{j-1} (\omega)) \leq \int_\Omega f^\star \text{ d}\overline{P}_{\tilde{\mathcal{E}}}.
\end{equation}
$\underline{P}_{\tilde{\mathcal{E}}} -$a.s. If furthermore $\underline{P}_{\tilde{\mathcal{E}}}$ is convex, then
\begin{equation}\label{ergo.result_sum}
\inf_{P_{\tilde{\mathcal{E}}}\in\mathcal{P}^{\text{co}}_{\tilde{\mathcal{E}}}}\sum_{\omega \in \Omega} f^\star(\omega) {P}_{\tilde{\mathcal{E}}}(\{\omega\}) \leq \lim\limits_{k \rightarrow \infty} \frac{1}{k} \sum_{j=1}^{k} f(T^{j-1} (\omega)) \leq \sup_{P_{\tilde{\mathcal{E}}}\in\mathcal{P}^{\text{co}}_{\tilde{\mathcal{E}}}}\sum_{\omega \in \Omega} f^\star(\omega) {P}_{\tilde{\mathcal{E}}}(\{\omega\})
\end{equation}
$\underline{P}_{\tilde{\mathcal{E}}} -$a.s.
\end{theorem}

%eq51=ergo.result
%eq98=lim.emp.avg
In appendix \ref{app_lemmas}, we give sufficient conditions for $\underline{P}_{\tilde{\mathcal{E}}}$ to be ($T$-)invariant, ergodic, and convex in Lemmas \ref{suff_cond_invariance}, \ref{suff_cond_inv_ergo}, and \ref{suff_cond_stricter_than_before}, respectively.

%Let us now explain what do we mean by ``$f^\star$ is a version of the conditional expectation of $f$ given $\mathcal{G}$'' in 

We now derive in Theorem \ref{subadd_ergo} a subadditive ergodic theorem for $\underline{P}_{\tilde{\mathcal{E}}}$, that can be seen as the DIPK version of \cite[Theorem 3]{cerreia} with an extra result, namely equation \eqref{extra_subadd}. As a consequence, in Corollary \ref{sharp_ergo1}, we find a sharpening of Theorem \ref{ergo1} when some additional assumptions are met; it can be seen as the DIPK version of \cite[Corollary 2]{cerreia}.\footnote{Notice that Corollary \ref{sharp_ergo1}.(2) is stronger than \cite[Corollary 2.(2)]{cerreia} since we state the equality between the lower expectation of $f$ and $f^\star$, and not just the equality between the Choquet integrals; given the convexity assumption on $\underline{P}_{\tilde{\mathcal{E}}}$, the equality between the Choquet integrals is implied.} 
%Before giving these results, we need to introduce some more concepts. A generic lower probability $\nu$ is Recall that a generic lower probability $\nu$ is \textit{convex} if $\nu(A \cup B) + \nu(A \cap B) \geq \nu(A) + \nu(B)$, for all $A,B \in \mathcal{F}$, and \textit{continuous at $\Omega$} if $\lim_{k \rightarrow \infty} \nu(A_k)=\nu(\Omega)$, when $A_k \uparrow \Omega$. In addition, $\nu$ is

We call a sequence $(S_k)$ of $\mathcal{F}$-measurable random variables \textit{superadditive} if $S_{k+\ell} \geq S_k +S_\ell \circ T^k$, for all $k$ and all $\ell$. It is \textit{subadditive} if the opposite inequality holds. It is {additive} if it is both super- and subadditive. A characterization of an additive sequence is the following: $(S_k)$ is additive if and only if there exists an $\mathcal{F}$-measurable functional $f$ on $\Omega$ such that
\begin{equation}\label{char_add}
S_k=\sum_{j=1}^k f \circ T^{j-1}, \quad \text{for all } k \in \mathbb{N}.
\end{equation} 
If we consider $(S_k)$ as in \eqref{char_add} and we take its absolute value, that is, if we consider $(|S_k|)$, we obtain a subadditive sequence. Notice also that if $f$ in our characterization belongs to $B(\Omega,\mathcal{F})$, we have that there exists $\lambda \in \mathbb{R}$ such that
\begin{equation}\label{bd_seq}
-\lambda k \leq S_k(\omega) \leq \lambda k, \quad \text{for all } k \in \mathbb{N},\omega \in \Omega.
\end{equation}
Similarly, $-\lambda k \leq |S_k(\omega)| \leq \lambda k$, for all $k \in \mathbb{N}$ and all $\omega \in \Omega$.

Subadditive stochastic processes (SSPs) were introduced for the first time in \cite{hammersley}, where the authors showed that they arise naturally in various contexts, but particularly in the study of random flows in lattices such as first-passage percolation. They developed an ergodic theory for these processes, that was then perfected in \cite{hammersley2,kingman,kingman2} and reconciled with Doob's martingale theory \cite{doob} in \cite{smeltzer}. In \cite{steele}, the author points out how SSPs provide an effective and unified approach to convergence theorems for the empirical discrepancy function, and he gives three noteworthy examples \cite[Section 4]{steele}. In \cite{kingman}, the author introduces also superadditive stochastic processes and points out how, since stochastic process $(S_k)$ is superadditive if and only if $(-S_k)$ is subadditive, any theorem about subadditive processes translates at once into a corresponding result about superadditive processes. Subsequent work specifically on superadditive processes include \cite{akcoglu,fong} where the authors study ratio ergodic theorems for superadditive processes, \cite{akcoglu2} where the authors study stochastic ergodic theorems for superadditive processes, and \cite{akcoglu3} where the authors study multiparameter generalizations of Kingman's pointwise theorem for continuous-time processes \cite{kingman2}.

%$(\Omega,\mathcal{F})$ is a standard measurable space, 

\begin{theorem}\label{subadd_ergo}
If $(S_k)$ is a super- or subadditive sequence satisfying \eqref{bd_seq} and $\underline{P}_{\tilde{\mathcal{E}}}$ is functionally invariant, then there is $f^\star \in B(\Omega,\mathcal{G})$ such that 
$$\lim_{k \rightarrow \infty} \frac{1}{k}S_k(\omega)=f^\star(\omega) \quad \underline{P}_{\tilde{\mathcal{E}}}-a.s.$$
In addition,
\begin{enumerate}
\item If $\underline{P}_{\tilde{\mathcal{E}}}$ is convex and strongly invariant, and $(S_k)$ is superadditive, then 
$$\inf_{P_{\tilde{\mathcal{E}}}\in\mathcal{P}^{\text{co}}_{\tilde{\mathcal{E}}}}\sum_{\omega \in \Omega} f^\star(\omega) {P}_{\tilde{\mathcal{E}}}(\{\omega\})=\sup_{k \in \mathbb{N}}\inf_{P_{\tilde{\mathcal{E}}}\in\mathcal{P}^{\text{co}}_{\tilde{\mathcal{E}}}} \frac{1}{k} \sum_{\omega \in \Omega}  S_k(\omega) {P}_{\tilde{\mathcal{E}}}(\{\omega\}).$$
\item If $\underline{P}_{\tilde{\mathcal{E}}}$ is convex and strongly invariant, and $(S_k)$ is subadditive, then 
$$\sup_{P_{\tilde{\mathcal{E}}}\in\mathcal{P}^{\text{co}}_{\tilde{\mathcal{E}}}}\sum_{\omega \in \Omega} f^\star(\omega) {P}_{\tilde{\mathcal{E}}}(\{\omega\})=\inf_{k \in \mathbb{N}}\sup_{P_{\tilde{\mathcal{E}}}\in\mathcal{P}^{\text{co}}_{\tilde{\mathcal{E}}}} \frac{1}{k} \sum_{\omega \in \Omega} S_k(\omega) {P}_{\tilde{\mathcal{E}}}(\{\omega\}).$$
\item If $\underline{P}_{\tilde{\mathcal{E}}}$ is ergodic and $(S_k)$ is either super- or subadditive, then 
$$\int_\Omega f^\star \text{ d} \underline{P}_{\tilde{\mathcal{E}}} \leq \lim_{k \rightarrow \infty} \frac{1}{k} S_k(\omega) \leq \int_\Omega f^\star \text{ d} \overline{P}_{\tilde{\mathcal{E}}}$$
$\underline{P}_{\tilde{\mathcal{E}}} -$a.s.
\item If $\underline{P}_{\tilde{\mathcal{E}}}$ is ergodic and convex, and $(S_k)$ is either super- or subadditive, then 
\begin{equation}\label{extra_subadd}
    \inf_{P_{\tilde{\mathcal{E}}}\in\mathcal{P}^{\text{co}}_{\tilde{\mathcal{E}}}}\sum_{\omega \in \Omega} f^\star(\omega) {P}_{\tilde{\mathcal{E}}}(\{\omega\}) \leq \lim_{k \rightarrow \infty} \frac{1}{k} S_k(\omega) \leq \sup_{P_{\tilde{\mathcal{E}}}\in\mathcal{P}^{\text{co}}_{\tilde{\mathcal{E}}}}\sum_{\omega \in \Omega} f^\star(\omega) {P}_{\tilde{\mathcal{E}}}(\{\omega\})
\end{equation}
$\underline{P}_{\tilde{\mathcal{E}}} -$a.s.
\end{enumerate}
\end{theorem}

In appendix \ref{app_lemmas}, we give sufficient conditions for $\underline{P}_{\tilde{\mathcal{E}}}$ to be strongly invariant and functionally invariant in Lemma \ref{suff_cond_stricter_than_before2}.
%$(\Omega,\mathcal{F})$ be a standard measurable space, and

\begin{corollary}\label{sharp_ergo1}
Let $\underline{P}_{\tilde{\mathcal{E}}}$ be convex and strongly invariant. Then, for all $f \in B(\Omega,\mathcal{F})$ there exists $f^\star \in B(\Omega,\mathcal{G})$ such that
\begin{equation*}
\lim_{k \rightarrow \infty} \frac{1}{k} \sum_{j=1}^k f(T^{j-1}(\omega))=f^\star(\omega) \quad \underline{P}_{\tilde{\mathcal{E}}}-a.s.
\end{equation*}
In addition, the following are true
\begin{enumerate}
\item For every $P \in \mathcal{I}$, $f^\star$ is a version of the conditional expectation of $f$ given $\mathcal{G}$.
%, that is, $$f^\star \in \left\lbrace{ \mathbb{E}_P \left(f \mid \mathcal{G} \right)}\right\rbrace_{P \in \mathcal{I}}.$$
\item $\inf_{P_{\tilde{\mathcal{E}}} \in \mathcal{P}^{\text{co}}_{\tilde{\mathcal{E}}}}\sum_{\omega \in \Omega} f^\star(\omega) {P}_{\tilde{\mathcal{E}}}(\{\omega\})=\inf_{P_{\tilde{\mathcal{E}}} \in \mathcal{P}^{\text{co}}_{\tilde{\mathcal{E}}}}\sum_{\omega \in \Omega} f(\omega) {P}_{\tilde{\mathcal{E}}}(\{\omega\})$.
\item If $\underline{P}_{\tilde{\mathcal{E}}}$ is also ergodic, then 
$$\inf_{P_{\tilde{\mathcal{E}}} \in \mathcal{P}^{\text{co}}_{\tilde{\mathcal{E}}}}\sum\limits_{\omega \in \Omega} f(\omega) {P}_{\tilde{\mathcal{E}}}(\{\omega\}) \leq \lim\limits_{k \rightarrow \infty} \frac{1}{k} \sum_{j=1}^{k} f(T^{j-1} (\omega)) \leq \sup_{P_{\tilde{\mathcal{E}}} \in \mathcal{P}^{\text{co}}_{\tilde{\mathcal{E}}}}\sum\limits_{\omega \in \Omega} f(\omega) {P}_{\tilde{\mathcal{E}}}(\{\omega\})$$
$\underline{P}_{\tilde{\mathcal{E}}} -$a.s.
\end{enumerate}
\end{corollary}
Corollary \ref{sharp_ergo1}.(1) means the following. Recall that in general, given a finite or countable probability space $(\Omega,\mathcal{F},P)$, an $\mathcal{F}$-measurable function $f$ on $\Omega$ having finite expectation, and a sub-sigma-algebra $\mathcal{G}$ of $\mathcal{F}$, $f$ is typically not $\mathcal{G}$-measurable, since $\mathcal{G} \subset \mathcal{F}$. However, a conditional expectation of $f$ given $\mathcal{G}$, denoted as $\mathbb{E}(f\mid \mathcal{G})$, can be defined as a $\mathcal{G}$-measurable function such that 
\begin{equation}\label{integral_eq}
    \sum_{\omega\in G} \mathbb{E}(f\mid \mathcal{G})(\omega) P(\{\omega\})=\sum_{\omega\in G} f(\omega) P(\{\omega\}),
\end{equation}
for all $G\in\mathcal{G}$; $\mathbb{E}(f\mid \mathcal{G})$ is in general not unique. Corollary \ref{sharp_ergo1}.(1) means that for all $f\in B(\Omega,\mathcal{F})$, we can find $f^\star \in B(\Omega,\mathcal{G})$ such that $f^\star = \mathbb{E}(f\mid \mathcal{G})$, and \eqref{integral_eq} holds for all $P\in\mathcal{I}$.\footnote{Notice that, in general, if $\Omega$ is uncountable, equation \eqref{integral_eq} becomes $\int_G \mathbb{E}(f\mid \mathcal{G}) \text{ d}P=\int_G f \text{ d}P$, for all $G\in\mathcal{G}$.}

\subsection{A strong law of large numbers}
A consequence of Theorem \ref{ergo1} is a strong law of large numbers. 
%for our belief updating procedure.
Before stating it, we need to introduce two notions. We first generalize the concept of a stationary stochastic process by allowing the underlying probability measure to be a lower probability. We then present the shift map, a classic idea in dynamics and ergodic theory.
%The two definitions are both from \cite{cerreia}.

%nonadditive, that is, by allowing it to be a Choquet capacity. 

Denote by $\mathbf{f}\equiv (f_k)_{k \in \mathbb{N}} \in B(\Omega,\mathcal{F})^{\mathbb{N}}$ a sequence of bounded and $\mathcal{F}$-measurable functionals on $\Omega$, and call $\mathcal{T}:=\bigcap_{\ell \in \mathbb{N}} \sigma(f_\ell,f_{\ell+1},\ldots)$ the tail sigma-algebra. Given a generic lower probability $\nu$ on $(\Omega,\mathcal{F})$, $\mathbf{f}$ is \textit{stationary} if and only if, for all $k \in \mathbb{N}$, all $\ell \in \mathbb{N}_0$, and  all Borel subset $B \subset \mathbb{R}^{\ell+1}$,
$$\nu \left( \left\lbrace{ \omega \in \Omega : \left( f_k(\omega),\ldots,f_{k+\ell}(\omega) \right) \in B }\right\rbrace \right)= \nu \left( \left\lbrace{ \omega \in \Omega : \left( f_{k+1}(\omega),\ldots,f_{k+\ell+1}(\omega) \right) \in B }\right\rbrace \right).$$

Now, denote by $(\mathbb{R}^\mathbb{N},\sigma(\mathcal{C}))$ the measurable space of sequences endowed with the sigma-algebra generated by the algebra of cylinders. Also denote by $s:\mathbb{R}^\mathbb{N} \rightarrow \mathbb{R}^\mathbb{N}$ the \textit{shift transformation}
$$s(x_1,x_2,x_3,\ldots)=(x_2,x_3,x_4,\ldots), \quad \text{for all } x \in \mathbb{R}^\mathbb{N} \text{.}$$

The sequence $\mathbf{f}$ induces a (natural) measurable map between $(\Omega,\mathcal{F})$ and measurable space $(\mathbb{R}^\mathbb{N},\sigma(\mathcal{C}))$ defined by
$$\omega \mapsto \mathbf{f}(\omega):=(f_1(\omega),\ldots,f_k(\omega),\ldots) \text{.}$$

Given any lower probability $\nu$ on $(\Omega,\mathcal{F})$, we can then define the map $\nu^\mathbf{f}:\sigma(\mathcal{C}) \rightarrow [0,1]$ as
$$C \mapsto \nu^\mathbf{f}(C):=\nu\left( \mathbf{f}^{-1}(C) \right)\text{.}$$

We say that $\mathbf{f}$ is {ergodic} if and only if $\nu^\mathbf{f}$ is ergodic with respect to the shift transformation. We are now ready for the strong law of large numbers. It can be seen as the DIPK version of \cite[Theorem 4]{cerreia}, but with a fewer assumption. Indeed, thanks to Lemma \ref{suff_cond_stricter_than_before2}.(i), we know that $\underline{P}_{\tilde{\mathcal{E}}}$ is always continuous at $\Omega$, so we do not need to explicitly require it in the statement of the theorem.
\begin{theorem}\label{slln}
Let $\underline{P}_{\tilde{\mathcal{E}}}$ be convex. If $\mathbf{f}=(f_n)_{n \in \mathbb{N}}$ is stationary and ergodic, then
\begin{equation}\label{slln_eq_imp}
    \inf_{P_{\tilde{\mathcal{E}}} \in \mathcal{P}^{\text{co}}_{\tilde{\mathcal{E}}}}\sum_{\omega \in \Omega} f_1(\omega) {P}_{\tilde{\mathcal{E}}}(\{\omega\}) \leq \lim\limits_{k \rightarrow \infty} \frac{1}{k} \sum\limits_{j=1}^k f_j(\omega) \leq \sup_{P_{\tilde{\mathcal{E}}} \in \mathcal{P}^{\text{co}}_{\tilde{\mathcal{E}}}}\sum_{\omega \in \Omega} f_1(\omega) {P}_{\tilde{\mathcal{E}}}(\{\omega\})
\end{equation}
$\underline{P}_{\tilde{\mathcal{E}}} - $a.s.
\end{theorem}

Notice that the assumption of stationarity gives us the fact that the limit for $k$ growing to infinity of $\frac{1}{k} \sum_{j=1}^k f_j(\omega)$ exists $\underline{P}_{\tilde{\mathcal{E}}} -$a.s. Then, to characterize this limit in terms of the infimum and the supremum of the expected value of $f_1$ with respect to the probability measures in $\mathcal{P}^{\text{co}}_{\tilde{\mathcal{E}}}$, we need the ergodicity assumption.

%lowest and highest taken

In the imprecise probabilities literature there are two excellent papers that deal with strong law of large numbers for lower previsions, \cite{de_cooman2} and \cite{cozman}. In the first one, the authors prove weak and strong laws of large numbers for coherent lower previsions under very weak assumptions.\footnote{Lower prevision $\underline{P}$ is coherent if given an arbitrary subset $\mathscr{K}$ of the set of gambles on $\Omega$, we have that $\sup[\sum_{j=1}^n (X_j-\underline{P}(X_j)) -m (X_0-\underline{P}(X_0))] \geq 0$, whenever $m,n\in\mathbb{Z}_+$ and $X_0,X_1,\ldots,X_n$ (not necessarily distinct) are in $\mathscr{K}$ \cite[Section 2.5.1]{walley}.} Their results make precise the following statement, taken from \cite[Section 1]{de_cooman2}: suppose an agent gives a lower prevision $m$ for bounded random variables $X_1,\ldots,X_N$ and they assess that they cannot learn from the past, in the sense that observations of variables $X_1,\ldots,X_{k-1}$ do not affect the lower prevision for $X_k$, $k\in\{2,\ldots,N\}$. Then, coherence implies that they should not bet on the event that the sample mean dominates lower prevision $m$ at rates that increase to $1$ as the number of observations $N$ increases to infinity. There are two main differences between \cite{de_cooman2} and the results in this section. First, in \cite{de_cooman2} the authors work with lower previsions instead of lower probabilities and do not require the state space to be at most countable, a more general approach. We operate with lower probabilities and with a finite or countable $\Omega$ because the version of DIPK updating in the case when $\Omega$ is uncountable and the agent expresses their beliefs via previsions has not yet been developed. Second, in \cite{de_cooman2} the authors assume that the agent cannot learn from past observations by making the lower prevision for $X_k$ not change (staying equal to $m$) after observing outcomes $X_1,\ldots,X_{k-1}$. We do not require this: for us, the agent has already updated their beliefs on the elements of $\mathcal{F}$ via DIPK, beliefs that are now encapsulated in set $\mathcal{P}_{\tilde{\mathcal{E}}}^{\text{co}} \subset \Delta(\Omega,\mathcal{F})$. They then study a phenomenon represented by sequence $\mathbf{f}$ that, if stationary and ergodic, is such that \eqref{slln_eq_imp} holds $\underline{P}_{\tilde{\mathcal{E}}} - $a.s. We also point out that since the strong law of large numbers in \cite{de_cooman2} can be seen as a generalization of that in \cite{maccheroni}, which is one of the building blocks of Theorem \ref{slln}, we believe that in the future the results in \cite[Section 5]{de_cooman2} will prove fundamental to derive a strong law of large numbers for the more general version of DIPK updating.

In \cite{cozman}, the author builds upon the results in \cite{de_cooman2}. He proves laws of large numbers that require a weaker irrelevance assumption (a way of assuming that the agent cannot learn from the past) than the one in \cite{de_cooman2}, and that do not require random variables $X_1,\ldots,X_N$ to be bounded. The same differences highlighted between our work and \cite{de_cooman2} hold for the present paper and \cite{cozman}. This latter is more general in that it does not require $\Omega$ to be at most countable, and does not uses merely lower probabilities, but rather lower prevision (of possibly unbounded random variables). Also, \cite{cozman} works with a very weak form of the irrelevance assumption called \textit{weak forward regular irrelevance} \cite[Section 3]{cozman}, while we do not require irrelevance for the same motivations explained before. In the future, the results in \cite{cozman} -- especially Theorem 4 -- will prove useful in deriving a strong law of large numbers for the generalization of DIPK updating.

\section{Conclusion}\label{concl}
In this paper we give an ergodic theory for the limit $\mathcal{P}^\text{co}_{\tilde{\mathcal{E}}}$ of the sequence $(\mathcal{P}^\text{co}_{\mathcal{E}_n})$ of successive dynamic imprecise probability kinematics updates of a set $\mathcal{P}^\text{co}_{\mathcal{E}_0}$ of probabilities representing the initial beliefs of an agent on the elements of $\mathcal{F}=2^\Omega$, $\Omega$ assumed at most countable. A consequence of this ergodic theory is a strong law of large numbers. We believe these results are pivotal for achieving a wider applicability of DIPK, for example via a generalization of the classical MCMC procedure that would instead accommodates DIPK updating; this will be subject to forthcoming studies. In the future, we also plan to find sufficient conditions that are easier to verify than the ones we have given for $\underline{P}_{\tilde{\mathcal{E}}}$ to be ($T$-)invariant, ergodic, convex, strictly invariant, and functionally invariant in Lemmas \ref{suff_cond_invariance}, \ref{suff_cond_inv_ergo}, \ref{suff_cond_stricter_than_before}, and \ref{suff_cond_stricter_than_before2}, respectively. In addition, when a generalized version of DIPK updating that deals with uncountable $\Omega$ and with the agent expressing their beliefs via lower previsions will be developed, we aim at showing that a version of the results presented in this article continues to hold.
%a version of the results presented in this article continues to hold for the generalization of DIPK -- to deal with uncountable $\Omega$ and with the agent expressing their beliefs via lower previsions.     if for the generalization of DIPK -- to deal with uncountable $\Omega$ and with the agent expressing their beliefs via lower previsions -- we are able to prove that sequence $(\mathcal{P}^\text{co}_{\mathcal{E}_n})$ converges, then .

\section*{Declaration of Competing Interest}
The authors declare that they have no known competing financial interests or personal relationships that could have appeared to influence the work reported in this paper.

%state space $\Omega$, 

\section*{Acknowledgements}
We would like to express our gratitude to Simone Cerreia-Vioglio for discussions on conditional expected values and subadditive sequences of random variables, Teddy Seidenfeld for discussion on the number of extrema of the core of a lower probability and for helping with the proof of Lemma \ref{core_conv}. We would also like to thank three anonymous referees and Roberto Corrao for their generous suggestions regarding content and presentation. Michele Caprio would like to acknowledge partial funding from NSF CCF-1934964 and ARO MURI W911NF2010080.
Sayan Mukherjee would like to acknowledge partial funding from HFSP
RGP005, NSF DMS 17-13012, NSF BCS 1552848, NSF DBI 1661386, NSF IIS
15-46331, NSF DMS 16-13261, and the Alexander von Humboldt Foundation. 
Sayan Mukherjee would also like to acknowledge the German Federal Ministry of Education and Research
within the project Competence Center for Scalable Data Analytics and Artificial Intelligence (ScaDS.AI) Dresden/Leipzig (BMBF 01IS18026B).

\appendix
\section{Technical lemmas}\label{app_lemmas}
We first present sufficient conditions for $\underline{P}_{\tilde{\mathcal{E}}}$ to be $T$-invariant, ergodic, and convex.
\begin{lemma}\label{suff_cond_invariance}
If there exists $T\in\mathbb{N}_0:=\mathbb{N}\cup\{0\}$ such that for all $t\geq T$ we can always find a collection $\{P^A_{\mathcal{E}_t}\}_{A\in\mathcal{F}} \subset \mathcal{P}^{\text{co}}_{\mathcal{E}_t}$ such that for $A^\prime\in\mathcal{F}$
\begin{itemize}
    \item $P^{A^\prime}_{\mathcal{E}_t}(A^\prime)=P^{A^\prime}_{\mathcal{E}_t}(T^{-1}(A^\prime))$,
    \item $P^{A^\prime}_{\mathcal{E}_t}(A^\prime)\leq P_{\mathcal{E}_t}(A^\prime)$, for all $P_{\mathcal{E}_t}\in\mathcal{P}_{\mathcal{E}_t}$,
    \item $P^{A^\prime}_{\mathcal{E}_t}(T^{-1}(A^\prime))\leq P_{\mathcal{E}_t}(T^{-1}(A^\prime))$, for all $P_{\mathcal{E}_t}\in\mathcal{P}_{\mathcal{E}_t}$,
\end{itemize}
then $\underline{P}_{\tilde{\mathcal{E}}}$ is $T$-invariant.
\end{lemma}
%Notice that we require collection $\{P^A_{\mathcal{E}_t}\}_{A\in\mathcal{F}}$ to be a subset of $\mathcal{P}_{\mathcal{E}_t}$ because by Theorem \ref{choq_rel} we have that $\underline{P}_{\mathcal{E}_t}(A)=\inf_{P\in\mathcal{P}_{\mathcal{E}_t}}P(A)=\inf_{P\in\mathcal{P}^{\text{co}}_{\mathcal{E}_t}}P(A)$, for all $A\in\mathcal{F}$.

\begin{lemma}\label{suff_cond_inv_ergo}
If there exist $T\in\mathbb{N}_0$ and  $P_{\mathcal{E}_T}^\prime\in\mathcal{P}^{\text{co}}_{\mathcal{E}_T}$ such that $P_{\mathcal{E}_T}^\prime(A)=0$ for all $A\in\mathcal{G}$, then $\underline{P}_{\tilde{\mathcal{E}}}$ is ergodic.
\end{lemma}
In this lemma we are tacitly assuming that $\Omega$ does not belong to $\mathcal{G}$, otherwise no such $P_{\mathcal{E}_T}^\prime$ exists. That is, if we are dealing with an operator $T$ for which  $\Omega$ is in $\mathcal{G}$, then Lemma \ref{suff_cond_inv_ergo} holds with $\mathcal{G}^\prime:=\mathcal{G}\setminus\{\Omega\}$ in place of $\mathcal{G}$. A result of Lemma \ref{suff_cond_inv_ergo} is that if the agent selects $\mathcal{P}=\{P_1,\ldots,P_k\}$ such that it contains an element $P^\prime$ that assigns probability $0$ to all the elements in $\mathcal{G}$ (or $\mathcal{G}^\prime$), then $P^\prime$ belongs to $\mathcal{P}^\text{co}_{\mathcal{E}_0}$. In turn this ensures that $\underline{P}_{\tilde{\mathcal{E}}}$ is ergodic. 

Given two generic sets $A,B\in\mathcal{F}$ and a generic lower probability $\nu$ such that $\nu(B)\neq 0$, define the \textit{geometric conditional lower probability} $\nu^G(A\mid B)$ as $\nu^G(A\mid B):=\frac{\nu(A\cap B)}{\nu(B)}$. Then, we have the following.

\begin{lemma}\label{suff_cond_stricter_than_before}
If either of the following hold, then $\underline{P}_{\tilde{\mathcal{E}}}$ is convex.
\begin{itemize}
    \item[(i)] There exists $T\in\mathbb{N}_0$ such that for all $t\geq T$ and all $A,B\in\mathcal{F}$ such that $A\subset B$, there exists $P^\prime_{\mathcal{E}_t} \in \mathcal{P}^{\text{co}}_{\mathcal{E}_t}$ such that
    $$P^\prime_{\mathcal{E}_t}(A\mid E) = \underline{P}_{\mathcal{E}_t}^G(A\mid E) \quad \text{and} \quad P^\prime_{\mathcal{E}_t}(B\mid E) = \underline{P}_{\mathcal{E}_t}^G(B\mid E), \quad \text{for all } E\in\mathcal{E}_{t+1}.$$
    \item[(ii)] There exists $T\in\mathbb{N}_0$ such that for all $t\geq T$ and all finite chains $(A_i)_{i=1}^n \subset \mathcal{F}$, there exists $P^\prime_{\mathcal{E}_t} \in \mathcal{P}^{\text{co}}_{\mathcal{E}_t}$ such that
    $$P^\prime_{\mathcal{E}_t}(A_i\mid E) = \underline{P}_{\mathcal{E}_t}^G(A_i\mid E), \quad \text{for all } i\in\{1,\ldots,n\}\text{, for all } E\in\mathcal{E}_{t+1}.$$
    \item[(iii)] There exists $T\in\mathbb{N}_0$ such that for all $t\geq T$ and all chains $(A_i)_{i\in I} \subset \mathcal{F}$, there exists $P^\prime_{\mathcal{E}_t} \in \mathcal{P}^{\text{co}}_{\mathcal{E}_t}$ such that
    $$P^\prime_{\mathcal{E}_t}(A_i\mid E) = \underline{P}_{\mathcal{E}_t}^G(A_i\mid E), \quad \text{for all } i\in I\text{, for all } E\in\mathcal{E}_{t+1}.$$
\end{itemize}
\end{lemma}

Next, we give sufficient conditions for $\underline{P}_{\tilde{\mathcal{E}}}$ to be strongly invariant and functionally invariant.

\begin{lemma}\label{suff_cond_stricter_than_before2}
The following are true
\begin{itemize}
    \item[(i)] $\underline{P}_{\tilde{\mathcal{E}}}$ is always continuous at $\Omega$.
    \item[(ii)] If there exists $T\in\mathbb{N}_0$ such that for all $t\geq T$, $\underline{P}_{\mathcal{E}_t}$ is strongly invariant, then $\underline{P}_{\tilde{\mathcal{E}}}$ is strongly invariant.
    \item[(iii)] If there exists $T\in\mathbb{N}_0$ such that for all $t\geq T$, $P_{\mathcal{E}_t}(A\cap E)=P_{\mathcal{E}_t}(T^{-1}(A)\cap E)$, for all $A\in\mathcal{F}$, all $E\in\mathcal{E}_{t+1}$, and all $P_{\mathcal{E}_t} \in \mathcal{P}^{\text{co}}_{\mathcal{E}_t}$, then $\underline{P}_{\tilde{\mathcal{E}}}$ is functionally invariant.
\end{itemize}
\end{lemma}

We go on by showing how to construct a subadditive sequence starting from a superadditive sequence satisfying \eqref{bd_seq}. This will be used in the proof of Theorem \ref{subadd_ergo}.

\begin{lemma}\label{lemma1}
Let $(S_k)$ be a superadditive sequence satisfying \eqref{bd_seq}, and suppose $\mathcal{P}^{\text{co}}_{\tilde{\mathcal{E}}} \subset \mathcal{I}$. Define then the sequence $(a_k) \in \mathbb{R}^\mathbb{N}$ as $a_k:=- \inf_{P_{\tilde{\mathcal{E}}}\in\mathcal{P}^{\text{co}}_{\tilde{\mathcal{E}}}} \sum_{\omega \in \Omega} S_k(\omega) P_{\tilde{\mathcal{E}}}(\{\omega\})$, for all $k \in \mathbb{N}$. Then, $(a_k)$ is subadditive, that is, $a_{k+\ell} \leq a_k+a_\ell$, for all $k,\ell \in \mathbb{N}$. If $(S_k)$ is subadditive, we reach the same result by defining $a_k:= \sup_{P_{\tilde{\mathcal{E}}}\in\mathcal{P}^{\text{co}}_{\tilde{\mathcal{E}}}} \sum_{\omega \in \Omega} S_k(\omega) P_{\tilde{\mathcal{E}}}(\{\omega\})$.
\end{lemma}

The following, which is a direct consequence of \cite[Lemma 1]{cerreia}, gives sufficient conditions for $\underline{P}_{\tilde{\mathcal{E}}}^{\mathbf{f}}$ to be convex, continuous at $\mathbb{R}^\mathbb{N}$, and shift invariant, and for $\mathbf{f}$ to be ergodic.

\begin{lemma}\label{lem_slln}
If $\underline{P}_{\tilde{\mathcal{E}}}$ is convex and $\mathbf{f}$ is stationary, then $\underline{P}_{\tilde{\mathcal{E}}}^{\mathbf{f}}$ is convex, continuous at $\mathbb{R}^\mathbb{N}$, and shift invariant. In addition, $\underline{P}_{\tilde{\mathcal{E}}}(\mathcal{T})=\{0,1\}$ implies that $\mathbf{f}$ is ergodic.
\end{lemma}

In the next lemma we show that assuming the number of extrema of the core to be finite will not be a strong requirement in the future, when DIPK will be generalized to work with lower previsions of (bounded) random variables. Consider a finite collection of probability measures $\{P_1,\ldots,P_k\}$. 
%and let $\underline{P}$ be the associated lower probability. all $B(\Omega,\mathcal{F})$ be the space of bounded random variables on $\Omega$ and c 
Denote by $P(X)=\sum_{\omega\in\Omega} X(\omega) P(\{\omega\})$, for all $P\in \Delta(\Omega,\mathcal{F})$ and all $X\in B(\Omega,\mathcal{F})$, $\underline{P}(X)=\inf_{P\in\text{Conv}(P_1,\ldots,P_k)} P(X)$, and $\overline{P}(X)=\sup_{P\in\text{Conv}(P_1,\ldots,P_k)} P(X)$. Then, call
\begin{align*}
    \text{core}_B(\underline{P})&:=\{P\in\Delta(\Omega,\mathcal{F}) : P(X) \geq \underline{P}(X) \text{, for all } X\in B(\Omega,\mathcal{F})\}\\
    &=\{P\in\Delta(\Omega,\mathcal{F}) : \overline{P}(X) \geq P(X) \geq \underline{P}(X) \text{, for all } X\in B(\Omega,\mathcal{F})\}.
\end{align*}
\begin{lemma}\label{core_conv}
    The following is true
    $$\text{core}_B(\underline{P})=\text{Conv}(P_1,\ldots,P_k).$$
\end{lemma}
We conjecture that this lemma holds even when $\Omega$ is uncountable; this will be studied in future work.

We conclude with a result that will be used in the proof of Theorem \ref{ergo.result}. Before stating it, let us introduce some concepts from \cite{marinacci_ambig}. Given a generic measurable space $(\Omega,\mathcal{F})$, we call \textit{game} a real-valued set function $\nu$ on $\mathcal{F}$ such that $\nu(\emptyset)=0$. It is a \textit{charge} if $\nu(A\cup B)=\nu(A)+\nu(B)$ for all pairwise disjoint sets $A,B\in\mathcal{F}$. We call $ba(\mathcal{F})$ the set of charges on $\mathcal{F}$ having finite \textit{total variation norm} $\|\nu\|_{TV}:=\sup\sum_{i=1}^n |\nu(A_i)-\nu(A_{i-1})|$, where the supremum is taken over all finite chains $\emptyset=A_0 \subset A_1 \subset \cdots \subset A_n=\Omega$. Given a generic game $\nu$,
\begin{itemize}
    \item its \textit{core} -- that we denote $\text{core}_g(\nu)$ to stress the fact that is the core of a game $\nu$ -- is the set of all suitably normalized charges that setwise dominate it, that is,
$$\text{core}_g(\nu):=\left\lbrace{ \mu\in ba(\mathcal{F}) : \mu(A) \geq \nu(A) \text{, for all } A\in\mathcal{F} \text{, and } \mu(\Omega)=\nu(\Omega)}\right\rbrace;$$
    \item it is \textit{balanced} if $\text{core}_g(\nu)\neq\emptyset$, and \textit{exact} if it is balanced and $\nu(A)=\inf_{\mu\in\text{core}_g(\nu)}\mu(A)$, for all $A\in\mathcal{F}$;
    \item it is \textit{continuous at $\Omega$} if $\lim_{n\rightarrow\infty} \nu(A_n) =\nu(\Omega)$ whenever $A_n \uparrow \Omega$; %sequence $(A_n) \subset \mathcal{F}^\mathbb{N}$ is such that
    \item it is \textit{continuous} if it is \textit{continuous at each $A$}, that is, if for all $A\in\mathcal{F}$, $\lim_{n\rightarrow\infty} \nu(A_n) =\nu(A)$ whenever $A_n\uparrow A$ or $A_n\downarrow A$.
\end{itemize}
As we can see, these notions are very general and encompass the ones we gave in section \ref{NAD_imp}. An immediate consequence of \cite[Theorem 10]{marinacci_ambig} and the unnumbered remark following it is given in the next claim.
\begin{claim}\label{claim_imp_marin}
If $\nu$ is a positive exact game on $\mathcal{F}$, then it is continuous at $\Omega$ if and only if it is continuous at each $A$.
\end{claim}
In turn, Claim \ref{claim_imp_marin} implies the following.
\begin{lemma}\label{lemma_marinacci_nice}
  $\underline{P}_{\tilde{\mathcal{E}}}$ is always continuous.
\end{lemma}

%In turn, Lemma \ref{lemma_marinacci_nice} implies that since by Lemma \ref{suff_cond_stricter_than_before2} $\underline{P}_{\tilde{\mathcal{E}}}$ is always continuous at $\Omega$, then it is also continuous at each $A$, that is, $\underline{P}_{\tilde{\mathcal{E}}}$ is continuous in the sense of section \ref{NAD_imp}.(iii).

 %Game $\nu$ is said to be  Finally, game $\nu$ is , and 
\section{Dynamic (precise) probability kinematics}\label{dpk}
In this section, we introduce the dynamic (precise) probability kinematics (DPK) updating procedure. Suppose that $P$ is a generic probability measure on $(\Omega,\mathcal{F})$ -- $\Omega$ at most countable -- representing an agent's initial beliefs around the elements of $\mathcal{F}=2^\Omega$, and that we want to update it after collecting some data. The agent observes data points $x_1,\ldots,x_n$ that are realizations of a random quantity $X:\Omega \rightarrow \mathcal{X}$ whose distribution is unknown. 
%\footnote{$X$ is a fixed function.} 
%Notice that collecting $x_1,\ldots,x_n$ is equivalent to observing $\omega_1,\ldots,\omega_n \sim Q$, where $Q$ is unknown, and then computing $X(\omega_i)=x_i$. 
Consider now the collection $\mathcal{E}^\prime:=\{E_i\}_{i=1}^n$, where $E_i\equiv X^{-1}(x_i):=\{\omega \in \Omega : X(\omega)=x_i\}$. It induces partition $\mathcal{E}=\{E_j\}_{j=1}^{m+1}$ of $\Omega$, $m\leq n$, whose first $m$ elements are the unique elements of $\mathcal{E}^\prime$, and $E_{m+1}=\Omega\setminus\cup_{j=1}^m E_j$. The DPK update of $P$ is 
\begin{align}
P_{\mathcal{E}}:\mathcal{F} \rightarrow [0,1] , \quad A \mapsto &P_{\mathcal{E}}(A) :=\sum_{\emptyset \neq E_j \in \mathcal{E}} P(A \mid E_j) P_{\mathcal{E}}(E_j) \label{eq_dipk_prima}\\ 
\text{such that }& P_{\mathcal{E}}(E_j) \geq 0 \text{, for all } E_j \in \mathcal{E} \text{, and } \sum_{E_j \in \mathcal{E}} P_{\mathcal{E}}(E_j)=1. \nonumber
\end{align}
In particular, for all $E_j \in \mathcal{E}$, $P_{\mathcal{E}}(E_j)$ is defined as follows. Consider the empirical probability measure $P^{emp}\in\Delta(\Omega,\mathcal{F})$ such that, if $E_{m+1}\neq\emptyset$,
\begin{align*}
    P^{emp}(E_j)=\frac{1}{n+1} \sum_{i=1}^{n} \mathbb{I}(E_j=E_i), \quad \text{for all } j \in \{1,\ldots,m\},
\end{align*}
where $\mathbb{I}$ denotes the indicator function, and 
\begin{align*}
   P^{emp}(E_{m+1})=1-\sum_{j=1}^{m} P^{emp}(E_j).
\end{align*}
If instead $E_{m+1}=\emptyset$,
\begin{align*}
    P^{emp}(E_j)=\frac{1}{n} \sum_{i=1}^{n} \mathbb{I}(E_j=E_i), \quad \text{for all } j \in \{1,\ldots,m\}
\end{align*}
and 
\begin{align*}
   P^{emp}(E_{m+1})=0.
\end{align*}
Then, $P_\mathcal{E}$ is such that
\begin{align*}
    P_{\mathcal{E}}(E_j)=\beta(n) P(E_j) + \left[ 1-\beta(n) \right] P^{emp}(E_j), \quad \text{for all } E_j\in\mathcal{E},
\end{align*}
where $\beta(n)$ is a coefficient in $[0,1]$ depending on $n$, chosen by the agent. It is easy to see that $P_\mathcal{E}$ satisfies the Kolmogorovian axioms of probability measures.
%: the posterior probability $P_\mathcal{E}$ assigned to the elements $E_j$ of partition $\mathcal{E}$ is a mixture of the prior $P$ and the empirical probability measure $P^{emp}$.

Let us give a remark. We tacitly assumed that for all nonempty $A\in\mathcal{F}$, the probability assigned to $A$ by $P$ (representing the agent's initial beliefs) is positive. In formulas, 
\begin{equation}\label{tacit_assumption}
    P(A)>0 , \quad \text{ for all } \emptyset\neq A \in\mathcal{F}.
\end{equation}
This assumption is not too stringent: suppose the agent specifies $P$ so that there is a collection of sets $\{A^\prime_k\}\subset \mathcal{F}$ such that $A^\prime_k\neq\emptyset$ and $P(A^\prime_k)=0$, for all $k$. Then, we prescribe the agent to modify slightly their initial beliefs as follows. They should pick an arbitrary small $\epsilon>0$ and a collection of (nonempty and not necessarily different) elements $\{A^{\prime\prime}_k\}\subset\mathcal{F}$ such that $P(A^{\prime\prime}_k)>0$ for all $k$, and define a new probability measure $\tilde{P}$ such that $\tilde{P}(A^\prime_k)=P(A^\prime_k)+\epsilon=\epsilon$ and $0 < \tilde{P}(A^{\prime\prime}_k)=P(A^{\prime\prime}_k)-\epsilon$, for all $k$. This procedure keeps the initial beliefs essentially unaltered, and avoids complications coming from conditioning on $0$ probability events. To see this, notice that \eqref{tacit_assumption} implies that $P(A\mid E_j)$ in \eqref{eq_dipk_prima} is well defined for all $\emptyset\neq E_j\in\mathcal{E}$. In addition, the same complications are avoided in successive updating. Call $\mathcal{E}_2$ the partition resulting from observing new data points, and notice that $\mathcal{E}_2$ is not coarser than $\mathcal{E}\equiv\mathcal{E}_1$. Then, for all $A\in\mathcal{F}$,
\begin{align*}
    P_{\mathcal{E}_2}(A)&=\sum_{\emptyset\neq E_j\in\mathcal{E}_2} P_{\mathcal{E}_1}(A\mid E_j) P_{\mathcal{E}_2}(E_j)\\
    &=\sum_{\emptyset\neq E_j\in\mathcal{E}_2} \frac{P_{\mathcal{E}_1}(A\cap E_j)}{P_{\mathcal{E}_1}(E_j)} P_{\mathcal{E}_2}(E_j).
\end{align*}
Because $\mathcal{E}_2$ is not coarser than $\mathcal{E}_1$, then for all $E^2\in\mathcal{E}_2$ we can always find $E^1\in\mathcal{E}_1$ such that $E^1\supset E^2$. This implies that 
$$P_{\mathcal{E}_1}(E^2)=\frac{P(E^2)}{P(E^1)}\left[ \beta(n) P(E^1) + \left(1-\beta(n)\right) P^{emp}(E^1) \right],$$
which is positive by \eqref{tacit_assumption}. In turn, this gives us that $P_{\mathcal{E}_1}(A\mid E_j)$ is a well defined operation, for all $E_j\in\mathcal{E}_2$. A similar argument holds for successive updating based on partitions $\mathcal{E}_t$, $t>2$. 
%A similar argument shows that our  assumption prevents conditioning on $0$ probability events in successive updates. 
In the future, we plan to deal with the delicate matter of conditioning on $0$ probability events in a more sophisticated way, possibly using techniques from the literature on lexicographic probabilities \cite{blume} or layers of zero probabilities \cite{coletti}.

Notice also that \eqref{tacit_assumption} implies a \textit{near-ignorance} assumption in the DIPK update. This means that every element in the collection $\{P_1,\ldots,P_k\}$ specified by the agent at the beginning of the analysis gives positive probability to all nonempty $A \in\mathcal{F}$. This is desirable because no finite sample is enough to annihilate a sufficiently extreme prior belief. To see this, suppose that there is a $P\in\mathcal{P}^{\text{co}}_{\mathcal{E}_0}$ and an $A^\prime\in \mathcal{F}$ such that $P(A^\prime)=0$; then 
\begin{itemize}
    \item $\underline{P}(A^\prime)=0$, and
    \item $P_\mathcal{E}(A^\prime)=0$ as well, since $P(A^\prime\cap E) \leq P(A^\prime)$, for all $E\in\mathcal{E}$, by the monotonicity of probability measures. This implies that $\underline{P}_\mathcal{E}(A^\prime)=0$.
\end{itemize}
As we can see, no finite amount of data can resolve vacuous initial beliefs.

%This means that no finite sample is enough to annihilate a sufficiently extreme prior belief. There is then a compromise to be made, and this is the compromise of \textit{near-ignorance}. The near-ignorance collection $\{P_1,\ldots,P_k\}$ specified by the agent at the beginning of the analysis can be vacuous a priori in some respects, but not on all the elements of $\mathcal{F}$. The agent should specify probability measures that attach non-zero probability to the events they care the most about.

\section{Proofs}\label{proofs}
%Recall that under the same regularity conditions of \cite[Proposition 12]{prob.kin}, the sequence $(\mathcal{P}_{\mathcal{E}_t})$ such that $\mathcal{P}_{\mathcal{E}_t}:=\text{ex}(\mathcal{P}^{\text{co}}_{\mathcal{E}_t})$, for all $t$, converges to $\mathcal{P}_{\tilde{\mathcal{E}}}:=\text{ex}(\mathcal{P}^{\text{co}}_{\tilde{\mathcal{E}}})$. In addition, by Theorem \ref{choq_rel}.(c), we have that the lower probabilities associated with $\mathcal{P}_{\mathcal{E}_t}$ and $\mathcal{P}^{\text{co}}_{\mathcal{E}_t}$ coincide, for all $t$.

\begin{proof}[Proof of Theorem \ref{ergo1}]
From \cite[Corollary 1]{cerreia}, we know that if $\underline{P}_{\tilde{\mathcal{E}}}$ is invariant, then core$(\underline{P}_{\tilde{\mathcal{E}}})=\mathcal{P}^{\text{co}}_{\tilde{\mathcal{E}}} \subset \mathcal{PI}$, where $\mathcal{PI}$ is the set of potentially invariant probability measures; that is, a probability measure $P$ belongs to $\mathcal{PI}$ if and only if
$$\exists \hat{P} \in \mathcal{I} : P(E)=\hat{P}(E), \quad \text{for all } E \in \mathcal{G}.$$
Then, \cite[Theorem 5]{cerreia} ensures us that $\mathcal{P}^{\text{co}}_{\tilde{\mathcal{E}}} \subset \mathcal{PI}$ is equivalent to the fact that for all $f \in B(\Omega,\mathcal{F})$, there exists $f^\star \in B(\Omega,\mathcal{G})$ such that
$$\lim\limits_{k \rightarrow \infty} \frac{1}{k} \sum\limits_{j=1}^{k} f(T^{j-1}(\omega))=f^\star(\omega) \quad \underline{P}_{\tilde{\mathcal{E}}}-a.s.$$
In particular, $f^\star \in B(\Omega,\mathcal{G})$ is defined as
$$\omega \mapsto f^\star(\omega):=\limsup_{k \rightarrow \infty} \frac{1}{k} \sum_{j=1}^k f(T^{j-1}(\omega)).$$
So we retrieve \eqref{lim.emp.avg}: the limit of the empirical averages exists and is finite $\underline{P}_{\tilde{\mathcal{E}}} -$a.s.

Assume that $\underline{P}_{\tilde{\mathcal{E}}}$ is also ergodic, and suppose for now that $f^\star \geq 0$. Since $\underline{P}_{\tilde{\mathcal{E}}}$ is a lower probability such that $\underline{P}_{\tilde{\mathcal{E}}}(\mathcal{G})= \{0,1\}$, and $0 \leq f^\star \leq \lambda$ for some $\lambda \in \mathbb{R}$ (because $f^\star$ is bounded), then
$$I:=\left\lbrace{ t \in \mathbb{R}_+ : \underline{P}_{\tilde{\mathcal{E}}} \left( \{ \omega \in \Omega : f^\star(\omega)\geq t \} \right) =1 }\right\rbrace$$
and 
$$J:=\left\lbrace{ t \in \mathbb{R}_- : \underline{P}_{\tilde{\mathcal{E}}} \left( \{ \omega \in \Omega : -f^\star(\omega)\geq t \} \right) =1 }\right\rbrace$$
are well defined nonempty intervals. $I$ is bounded from above and such that $0 \in I$, and $J$ is bounded from below and such that $-\lambda \in J$. Notice also that by Lemma  \ref{lemma_marinacci_nice} $\underline{P}_{\tilde{\mathcal{E}}}$ is continuous. 
%since $\underline{P}_{\tilde{\mathcal{E}}}$ is a lower probability, then it is continuous. 
We can conclude that $\sup I=:t^\star \in I$ and $\sup J =:t_\star \in J$. Since $\underline{P}_{\tilde{\mathcal{E}}}(\mathcal{G}) = \{0,1\}$, we have that
\begin{align*}
\int_\Omega f^\star \text{ d} \underline{P}_{\tilde{\mathcal{E}}} = \int_0^\infty \underline{P}_{\tilde{\mathcal{E}}} \left( \left\lbrace{ \omega\in\Omega : f^\star( \omega) \geq t }\right\rbrace \right) \text{ d}t = \int_0^{\sup I} \text{d}t= t^\star,
\end{align*}
and 
\begin{align*}
\int_\Omega (-f^\star) \text{ d} \underline{P}_{\tilde{\mathcal{E}}}= \int_{-\infty}^0 \left[\underline{P}_{\tilde{\mathcal{E}}} \left( \left\lbrace{ \omega\in\Omega : -f^\star( \omega) \geq t }\right\rbrace \right) - \underline{P}_{\tilde{\mathcal{E}}}(\Omega) \right] \text{ d}t = \int_{\sup J}^0 (-1)\text{ d}t= t_\star.
\end{align*}
So we have that  $t^\star=\int_\Omega f^\star \text{ d} \underline{P}_{\tilde{\mathcal{E}}}$ and $t_\star=\int_\Omega (-f^\star) \text{ d} \underline{P}_{\tilde{\mathcal{E}}}$. Now, since $t^\star \in I$ and $t_\star \in J$, we also have that
$$\underline{P}_{\tilde{\mathcal{E}}} \left( \{\omega \in \Omega : f^\star( \omega) \geq t^\star\}\right)=1=\underline{P}_{\tilde{\mathcal{E}}} \left( \{\omega \in \Omega : f^\star( \omega) \leq -t_\star\}\right) \text{.}$$
Since $\underline{P}_{\tilde{\mathcal{E}}}$ is a lower probability, this implies that
\begin{align}\label{eq53}
\begin{split}
\underline{P}_{\tilde{\mathcal{E}}} \left( \left\lbrace{ \omega \in \Omega : \int_\Omega f^\star \text{ d} \underline{P}_{\tilde{\mathcal{E}}} \leq f^\star (\omega) \leq \int_\Omega f^\star \text{ d} \overline{P}_{\tilde{\mathcal{E}}} }\right\rbrace\right)=\underline{P}_{\tilde{\mathcal{E}}} \left( \left\lbrace{ \omega \in \Omega : t^\star \leq f^\star(\omega) \leq -t_\star }\right\rbrace\right)=1 \text{.}
\end{split}
\end{align}
Let us now relax the assumption that $f^\star \geq 0$. Since $f^\star \in B(\Omega,\mathcal{G})$, then there exists $c \in \mathbb{R}$ such that $f^\star+c\mathbbm{1}_\Omega \geq 0$. By \eqref{eq53}, we have that
\begin{equation*}
\begin{split}
\underline{P}_{\tilde{\mathcal{E}}} &\bigg( \bigg\{ \omega \in \Omega :
\int_\Omega (f^\star+c\mathbbm{1}_\Omega) \text{ d} \underline{P}_{\tilde{\mathcal{E}}} \leq f^\star (\omega)+c \leq \int_\Omega (f^\star+c\mathbbm{1}_\Omega) \text{ d} \overline{P}_{\tilde{\mathcal{E}}}\bigg\}\bigg)\\
=\underline{P}_{\tilde{\mathcal{E}}} &\bigg( \bigg\{ \omega \in \Omega :
\int_\Omega f^\star \text{ d} \underline{P}_{\tilde{\mathcal{E}}} +c \leq f^\star (\omega)+c \leq \int_\Omega f^\star \text{ d} \overline{P}_{\tilde{\mathcal{E}}} +c\bigg\}\bigg)\\
=\underline{P}_{\tilde{\mathcal{E}}} &\bigg( \bigg\{ \omega \in \Omega :
\int_\Omega f^\star \text{ d} \underline{P}_{\tilde{\mathcal{E}}}  \leq f^\star (\omega) \leq \int_\Omega f^\star \text{ d} \overline{P}_{\tilde{\mathcal{E}}} \bigg\}\bigg)=1.
\end{split}
\end{equation*}
To conclude the proof, since by \eqref{lim.emp.avg} we have that 
$$\underline{P}_{\tilde{\mathcal{E}}} \left( \left\lbrace{\omega \in \Omega : f^\star(\omega)=\lim_{k \rightarrow \infty} \sum_{j=1}^k f(T^{j-1}(\omega))}\right\rbrace \right)=1$$
and since $\underline{P}_{\tilde{\mathcal{E}}}$ is a lower probability, this implies that 
\begin{align}\label{equaz_imp1}
    \underline{P}_{\tilde{\mathcal{E}}} &\bigg( \bigg\{ \omega \in \Omega :
\int_\Omega f^\star \text{ d} \underline{P}_{\tilde{\mathcal{E}}}  \leq \lim_{k \rightarrow \infty} \sum_{j=1}^k f(T^{j-1}(\omega)) \leq \int_\Omega f^\star \text{ d} \overline{P}_{\tilde{\mathcal{E}}} \bigg\}\bigg)=1,
\end{align}
retrieving equation \eqref{ergo.result_int}. If furthermore $\underline{P}_{\tilde{\mathcal{E}}}$ is convex, we have that by \cite[Theorem 38]{marinacci_ambig}
\begin{equation}\label{equaz_imp2}
    \int_\Omega f^\star \text{ d} \underline{P}_{\tilde{\mathcal{E}}}=\inf_{P_{\tilde{\mathcal{E}}}\in\mathcal{P}^{\text{co}}_{\tilde{\mathcal{E}}}}\int_\Omega f^\star \text{ d} {P}_{\tilde{\mathcal{E}}}
\end{equation}
and that, since $\Omega$ is at most countable, 
\begin{equation}\label{equaz_imp3}
    \int_\Omega f^\star \text{ d} {P}_{\tilde{\mathcal{E}}}=\sum_{\omega\in\Omega} f^\star(\omega) {P}_{\tilde{\mathcal{E}}}(\{\omega\}).
\end{equation}
Now, substituting \eqref{equaz_imp3} in \eqref{equaz_imp2}, we obtain 
\begin{equation}\label{equaz_imp4}
    \int_\Omega f^\star \text{ d} \underline{P}_{\tilde{\mathcal{E}}}=\inf_{P_{\tilde{\mathcal{E}}}\in\mathcal{P}^{\text{co}}_{\tilde{\mathcal{E}}}}\sum_{\omega\in\Omega} f^\star(\omega) {P}_{\tilde{\mathcal{E}}}(\{\omega\}).
\end{equation}
Also, as a consequence, $\int_\Omega f^\star \text{ d} \overline{P}_{\tilde{\mathcal{E}}}=\sup_{P_{\tilde{\mathcal{E}}}\in\mathcal{P}^{\text{co}}_{\tilde{\mathcal{E}}}}\sum_{\omega\in\Omega} f^\star(\omega) {P}_{\tilde{\mathcal{E}}}(\{\omega\})$. Plugging this and \eqref{equaz_imp4} in \eqref{equaz_imp1}, we obtain \eqref{ergo.result_sum}, concluding the proof.
\end{proof}

\begin{proof}[Proof of Theorem \ref{subadd_ergo}]
Given our assumption that $\underline{P}_{\tilde{\mathcal{E}}}$ is functionally invariant, we have that by \cite[Theorem 1]{cerreia} $\mathcal{P}^{\text{co}}_{\tilde{\mathcal{E}}} \subset \mathcal{I}$. Define the sequence $(f_k)$ such that $f_k(\omega):=\frac{1}{k} S_k(\omega)$, for all $k \in \mathbb{N}$ and all $\omega \in \Omega$. By $(S_k)$ satisfying \eqref{bd_seq}, this implies that $f_k \in B(\Omega,\mathcal{F})$, for all $k$. Consider then a function $p:\mathcal{F}\times\Omega \rightarrow [0,1]$ such that
\begin{itemize}
\item for all $P_{\tilde{\mathcal{E}}} \in \mathcal{P}^{\text{co}}_{\tilde{\mathcal{E}}}$ and all $A \in \mathcal{F}$, $p(A,\cdot):\Omega \rightarrow [0,1]$ is a version of the conditional probability of $A$ given $\mathcal{G}$;
%that is, $p(A,\cdot) \in P(A \mid \mathcal{G})$;
\item for all $\omega \in \Omega$, $p(\cdot,\omega):\mathcal{F} \rightarrow [0,1]$ is a probability measure;
\item for all $\omega \in \Omega$, $p(\cdot,\omega) \in \mathcal{P}^{\text{co}}_{\tilde{\mathcal{E}}}$.
\end{itemize}
For all $k \in \mathbb{N}$, define then 
\begin{equation}\label{f_hat}
\hat{f}_k:\Omega \rightarrow \mathbb{R}, \quad \omega \mapsto \hat{f}_k(\omega):=\sum_{\omega^\prime \in \Omega} f_k(\omega^\prime) p(\{\omega^\prime\},\omega).
\end{equation}
%Notice that, for all $k$, $\hat{f}_k \in \cap_{P \in \mathcal{P}_{\tilde{\mathcal{E}}}} \mathbb{E}_P(f_k \mid \mathcal{G})$. 
Notice that because $\Omega$ is at most countable, we can write $\hat{f}_k(\omega)$ as $\int_\Omega f \text{ d}p(\cdot,\omega)$, for all $\omega\in\Omega$. Given that $f_k \in B(\Omega,\mathcal{F})$, for all $k$, this implies that $\hat{f}_k \in B(\Omega,\mathcal{G})$, for all $k$. Since $(S_k)$ satisfies \eqref{bd_seq}, it follows that there exists $\lambda\in\mathbb{R}$ such that $-\lambda \leq f_k,\hat{f}_k \leq \lambda$, for all $k \in \mathbb{N}$. Define now $f^\star \in B(\Omega,\mathcal{G})$ by $f^\star:=\sup_{k \in \mathbb{N}} \hat{f}_k$. By Kingman's Subadditive Ergodic Theorem \cite[Theorem 10.7.1]{dudley} and \cite[Theorem 8.4]{gray}, we have that $f^\star=\lim_{k \rightarrow \infty} \hat{f}_k$ and $\lim_{k \rightarrow \infty} \frac{1}{k}S_k(\omega)=f^\star(\omega)$ $P_{\tilde{\mathcal{E}}} -$a.s., for all $P_{\tilde{\mathcal{E}}} \in \mathcal{P}^{\text{co}}_{\tilde{\mathcal{E}}}$. Since $\underline{P}_{\tilde{\mathcal{E}}}$ is a lower probability, it follows that $\lim_{k \rightarrow \infty} \frac{1}{k}S_k(\omega)=f^\star(\omega)$ $\underline{P}_{\tilde{\mathcal{E}}} -$a.s. This shows the first part of the theorem. Let us now show claims (1)--(4).

(1). If $\underline{P}_{\tilde{\mathcal{E}}}$ is strongly invariant, by \cite[Theorem 1]{cerreia} we have that core$(\underline{P}_{\tilde{\mathcal{E}}})=\mathcal{P}^{\text{co}}_{\tilde{\mathcal{E}}} \subset \mathcal{I}$. If it is also convex, by \cite[Theorem 38]{marinacci_ambig} we have that
\begin{equation}\label{imp_eq}
\inf_{P_{\tilde{\mathcal{E}}} \in \mathcal{P}^{\text{co}}_{\tilde{\mathcal{E}}}}\sum_{\omega \in \Omega} f(\omega) {P}_{\tilde{\mathcal{E}}}(\{\omega\})= \inf_{P_{\tilde{\mathcal{E}}} \in \mathcal{P}^{\text{co}}_{\tilde{\mathcal{E}}}} \int_{\Omega} f \text{ d}P_{\tilde{\mathcal{E}}} =  \int_{\Omega} f \text{ d}\underline{P}_{\tilde{\mathcal{E}}}, \quad \text{for all } f \in B(\Omega,\mathcal{F}),
\end{equation}
where the first equality comes from $\Omega$ being at most countable. Consider now the sequence $(a_k)$ defined as $a_k:=-\inf_{P_{\tilde{\mathcal{E}}} \in \mathcal{P}^{\text{co}}_{\tilde{\mathcal{E}}}} \sum_{\omega \in \Omega} S_k(\omega) {P}_{\tilde{\mathcal{E}}}(\{\omega\})$, for all $k \in \mathbb{N}$. By \eqref{imp_eq}, we have that  
$$a_k:=-\inf_{P_{\tilde{\mathcal{E}}} \in \mathcal{P}^{\text{co}}_{\tilde{\mathcal{E}}}} \sum_{\omega \in \Omega} S_k(\omega) {P}_{\tilde{\mathcal{E}}}(\{\omega\})=-\inf_{P_{\tilde{\mathcal{E}}} \in \mathcal{P}^{\text{co}}_{\tilde{\mathcal{E}}}} \int_{\Omega} S_k \text{ d}{P}_{\tilde{\mathcal{E}}}=- \int_{\Omega} S_k \text{ d}\underline{P}_{\tilde{\mathcal{E}}}, \quad \text{for all } k \in \mathbb{N}.$$
By Lemma \ref{lemma1}, it follows that $(a_k)$ is subadditive. By \cite[Lemma 8.3]{gray}, this implies that 
\begin{equation}\label{imp_eq2}
\lim_{k \rightarrow \infty} \frac{1}{k}(-a_k)=\sup_{k \in \mathbb{N}} \frac{1}{k}(-a_k).
\end{equation}
Now, by the fact that $(\hat{f}_n)$ is uniformly bounded, \cite[Theorem 22]{cerreia4}, equation \eqref{imp_eq2}, the first part of the theorem, and the fact that $\mathcal{P}^{\text{co}}_{\tilde{\mathcal{E}}} \subset \mathcal{I}$, the following equalities hold
\begin{align*}
\inf_{P_{\tilde{\mathcal{E}}} \in \mathcal{P}^{\text{co}}_{\tilde{\mathcal{E}}}}\sum_{\omega \in \Omega} f^\star(\omega) {P}_{\tilde{\mathcal{E}}}(\{\omega\}) &= \int_\Omega f^\star \text{ d} \underline{P}_{\tilde{\mathcal{E}}}=\int_{\Omega} \lim_{k \rightarrow \infty} \hat{f}_k \text{ d} \underline{P}_{\tilde{\mathcal{E}}} = \lim_{k \rightarrow \infty} \int_{\Omega} \hat{f}_k \text{ d} \underline{P}_{\tilde{\mathcal{E}}}\\
&= \lim_{k \rightarrow \infty} \left[\inf_{P_{\tilde{\mathcal{E}}} \in \mathcal{P}^{\text{co}}_{\tilde{\mathcal{E}}}}\int_{\Omega}  \hat{f}_k \text{ d} P_{\tilde{\mathcal{E}}} \right]=\lim_{k \rightarrow \infty} \left[\inf_{P_{\tilde{\mathcal{E}}} \in \mathcal{P}^{\text{co}}_{\tilde{\mathcal{E}}}}\int_{\Omega}  {f}_k \text{ d} P_{\tilde{\mathcal{E}}} \right]\\
&= \lim_{k \rightarrow \infty} \int_{\Omega} {f}_k \text{ d}\underline{P}_{\tilde{\mathcal{E}}}= \lim_{k \rightarrow \infty} \frac{1}{k} \int_{\Omega} S_k \text{ d}\underline{P}_{\tilde{\mathcal{E}}}\\
&= \lim_{k \rightarrow \infty} \frac{1}{k}(-a_k)=\sup_{k \in \mathbb{N}} \frac{1}{k}(-a_k)\\
&= \sup_{k \in \mathbb{N}}\frac{1}{k} \int_{\Omega} S_k \text{ d} \underline{P}_{\tilde{\mathcal{E}}}=\sup_{k \in \mathbb{N}} \inf_{P_{\tilde{\mathcal{E}}} \in \mathcal{P}^{\text{co}}_{\tilde{\mathcal{E}}}}\frac{1}{k}\sum_{\omega\in\Omega} S_k(\omega) {P}_{\tilde{\mathcal{E}}}(\{\omega\})
%&= \sup_{k \in \mathbb{N}} \sum_{\omega \in \Omega} {f}_k(\omega) \underline{P}_{\tilde{\mathcal{E}}}(\{\omega\}),
\end{align*}
concluding the proof of (1).

(2). If $\underline{P}_{\tilde{\mathcal{E}}}$ is strongly invariant, by \cite[Theorem 1]{cerreia} we have that core$(\underline{P}_{\tilde{\mathcal{E}}})=\mathcal{P}^{\text{co}}_{\tilde{\mathcal{E}}} \subset \mathcal{I}$. If it is also convex, by \cite[Theorem 38]{marinacci_ambig} we have that
\begin{equation}\label{imp_eq3}
\sup_{P_{\tilde{\mathcal{E}}} \in \mathcal{P}^{\text{co}}_{\tilde{\mathcal{E}}}}\sum_{\omega \in \Omega} f(\omega) {P}_{\tilde{\mathcal{E}}}(\{\omega\})= \sup_{P_{\tilde{\mathcal{E}}} \in \mathcal{P}^{\text{co}}_{\tilde{\mathcal{E}}}} \int_{\Omega} f \text{ d}P =  \int_{\Omega} f \text{ d}\overline{P}_{\tilde{\mathcal{E}}}, \quad \text{for all } f \in B(\Omega,\mathcal{F}).
\end{equation}
where the first equality comes from $\Omega$ being at most countable. Consider now the sequence $(a_k)$ defined as $a_k:=\sup_{P_{\tilde{\mathcal{E}}} \in \mathcal{P}^{\text{co}}_{\tilde{\mathcal{E}}}} \sum_{\omega \in \Omega} S_k(\omega) {P}_{\tilde{\mathcal{E}}}(\{\omega\})$, for all $k \in \mathbb{N}$. By \eqref{imp_eq3}, we have that  
$$a_k:=\sup_{P_{\tilde{\mathcal{E}}} \in \mathcal{P}^{\text{co}}_{\tilde{\mathcal{E}}}} \sum_{\omega \in \Omega} S_k(\omega) {P}_{\tilde{\mathcal{E}}}(\{\omega\})=\sup_{P_{\tilde{\mathcal{E}}} \in \mathcal{P}^{\text{co}}_{\tilde{\mathcal{E}}}} \int_{\Omega} S_k \text{ d}{P}_{\tilde{\mathcal{E}}}=\int_{\Omega} S_k \text{ d}\overline{P}_{\tilde{\mathcal{E}}}, \quad \text{for all } k \in \mathbb{N}.$$
By Lemma \ref{lemma1}, it follows that $(a_k)$ is subadditive. By \cite[Lemma 8.3]{gray}, this implies that 
\begin{equation}\label{imp_eq4}
\lim_{k \rightarrow \infty} \frac{a_k}{k}=\inf_{k \in \mathbb{N}} \frac{a_k}{k}.
\end{equation}
Now, by the fact that $(\hat{f}_n)$ is uniformly bounded, \cite[Theorem 22]{cerreia4}, equation \eqref{imp_eq4}, the first part of the theorem, and the fact that $\mathcal{P}^{\text{co}}_{\tilde{\mathcal{E}}} \subset \mathcal{I}$, the following equalities hold
\begin{align*}
\sup_{P_{\tilde{\mathcal{E}}} \in \mathcal{P}^{\text{co}}_{\tilde{\mathcal{E}}}}\sum_{\omega \in \Omega} f^\star(\omega) {P}_{\tilde{\mathcal{E}}}(\{\omega\}) &= \int_\Omega f^\star \text{ d} \overline{P}_{\tilde{\mathcal{E}}}=\int_{\Omega} \lim_{k \rightarrow \infty} \hat{f}_k \text{ d} \overline{P}_{\tilde{\mathcal{E}}} = \lim_{k \rightarrow \infty} \int_{\Omega} \hat{f}_k \text{ d} \overline{P}_{\tilde{\mathcal{E}}}\\
&= \lim_{k \rightarrow \infty} \left[\sup_{P_{\tilde{\mathcal{E}}} \in \mathcal{P}^{\text{co}}_{\tilde{\mathcal{E}}}}\int_{\Omega}  \hat{f}_k \text{ d} P_{\tilde{\mathcal{E}}} \right]=\lim_{k \rightarrow \infty} \left[\sup_{P_{\tilde{\mathcal{E}}} \in \mathcal{P}^{\text{co}}_{\tilde{\mathcal{E}}}}\int_{\Omega}  {f}_k \text{ d} P_{\tilde{\mathcal{E}}} \right]\\
&= \lim_{k \rightarrow \infty} \int_{\Omega} {f}_k \text{ d}\overline{P}_{\tilde{\mathcal{E}}}= \lim_{k \rightarrow \infty} \frac{1}{k} \int_{\Omega} S_k \text{ d}\overline{P}_{\tilde{\mathcal{E}}}\\
&= \lim_{k \rightarrow \infty} \frac{a_k}{k}=\inf_{k \in \mathbb{N}} \frac{a_k}{k}\\
&= \inf_{k \in \mathbb{N}}\frac{1}{k} \int_{\Omega} S_k \text{ d} \overline{P}_{\tilde{\mathcal{E}}}=\inf_{k \in \mathbb{N}} \sup_{P_{\tilde{\mathcal{E}}} \in \mathcal{P}^{\text{co}}_{\tilde{\mathcal{E}}}}\frac{1}{k}\sum_{\omega\in\Omega} S_k(\omega) {P}_{\tilde{\mathcal{E}}}(\{\omega\})
%&= \sup_{k \in \mathbb{N}} \sum_{\omega \in \Omega} {f}_k(\omega) \underline{P}_{\tilde{\mathcal{E}}}(\{\omega\}),
\end{align*}
concluding the proof of (2).

(3). If $\underline{P}_{\tilde{\mathcal{E}}}$ is ergodic, using the same technique as in the proof of Theorem \ref{ergo1} we can show that 
$$\underline{P}_{\tilde{\mathcal{E}}} \left( \left\lbrace{\omega \in \Omega : \int_{\Omega} f^\star \text{ d} \underline{P}_{\tilde{\mathcal{E}}}\leq f^\star(\omega) \leq \int_{\Omega} f^\star \text{ d} \overline{P}_{\tilde{\mathcal{E}}} }\right\rbrace\right)=1.$$
By the first part of the theorem, we have that 
$$\underline{P}_{\tilde{\mathcal{E}}} \left( \left\lbrace{\omega \in \Omega : \lim_{k \rightarrow \infty} \frac{1}{k}S_k(\omega) = f^\star(\omega) }\right\rbrace\right)=1.$$
Since $\underline{P}_{\tilde{\mathcal{E}}}$ is a lower probability, this implies that
\begin{equation}\label{equaz_imp5}
    \underline{P}_{\tilde{\mathcal{E}}} \left( \left\lbrace{\omega \in \Omega : \int_{\Omega} f^\star \text{ d} \underline{P}_{\tilde{\mathcal{E}}}\leq \lim_{k \rightarrow \infty} \frac{1}{k}S_k(\omega) \leq \int_{\Omega} f^\star \text{ d} \overline{P}_{\tilde{\mathcal{E}}} }\right\rbrace\right)=1.
\end{equation}

(4). If $\underline{P}_{\tilde{\mathcal{E}}}$ is ergodic and convex, we know that \eqref{equaz_imp4} and its consequence $\int_\Omega f^\star \text{ d} \overline{P}_{\tilde{\mathcal{E}}}=\sup_{P_{\tilde{\mathcal{E}}}\in\mathcal{P}^{\text{co}}_{\tilde{\mathcal{E}}}}\sum_{\omega\in\Omega} f^\star(\omega) {P}_{\tilde{\mathcal{E}}}(\{\omega\})$ hold. Plugging these in \eqref{equaz_imp5}, we obtain the fourth statement of the theorem, concluding the proof.
\end{proof}

\begin{proof}[Proof of Corollary \ref{sharp_ergo1}]
Consider $f \in B(\Omega,\mathcal{F})$. We first notice that $(S_k)$ such that $S_k:=\sum_{j=1}^k f \circ T^{j-1}$, for all $k \in \mathbb{N}$, is an additive sequence satisfying \eqref{bd_seq}. Now, since $\underline{P}_{\tilde{\mathcal{E}}}$ is 
%convex, continuous at $\Omega$, and 
strongly invariant, then it is also functionally invariant by \cite[Theorem 1]{cerreia}. Consider now the sequence $(f_k)$ where $f_k:=S_k/k$, for all $k \in \mathbb{N}$. Notice that $\hat{f}_k=\hat{f}$, for all $k \in \mathbb{N}$, where $\hat{f}_k$ is defined as in \eqref{f_hat}, $\hat{f}:\Omega \rightarrow \mathbb{R}$, $\omega \mapsto \hat{f}(\omega):=\sum_{\omega^\prime \in \Omega} f(\omega^\prime)p(\{\omega^\prime\},\omega)$, and $p:\mathcal{F}\times\Omega \rightarrow [0,1]$ is defined as in the proof of Theorem \ref{subadd_ergo}. Then, by the proof of Theorem \ref{subadd_ergo}, it follows that $\lim_{k \rightarrow \infty} \frac{1}{k} S_k = \lim_{k \rightarrow \infty} \hat{f}_k=\hat{f}$, $\underline{P}_{\tilde{\mathcal{E}}} -$a.s. This proves the first part of the corollary, and also point (1) by setting $f^\star=\hat{f}$. Let us now show claims (2) and (3).

(2). If $\underline{P}_{\tilde{\mathcal{E}}}$ is strongly invariant, by \cite[Theorem 1]{cerreia} we have that core$(\underline{P}_{\tilde{\mathcal{E}}})=\mathcal{P}^{\text{co}}_{\tilde{\mathcal{E}}} \subset \mathcal{I}$. If it is also convex, by \cite[Theorem 38]{marinacci_ambig} we have that
\begin{equation*}
\inf_{P_{\tilde{\mathcal{E}}} \in \mathcal{P}^{\text{co}}_{\tilde{\mathcal{E}}}}\sum_{\omega \in \Omega} f(\omega) {P}_{\tilde{\mathcal{E}}}(\{\omega\})= \inf_{P_{\tilde{\mathcal{E}}} \in \mathcal{P}^{\text{co}}_{\tilde{\mathcal{E}}}} \int_{\Omega} f \text{ d}P =  \int_{\Omega} f \text{ d}\underline{P}_{\tilde{\mathcal{E}}}, \quad \text{for all } f \in B(\Omega,\mathcal{F}),
\end{equation*}
where the first equality comes from $\Omega$ being at most countable. By (1) and $\mathcal{P}^{\text{co}}_{\tilde{\mathcal{E}}} \subset \mathcal{I}$, we have that 
\begin{align*}
\inf_{P_{\tilde{\mathcal{E}}} \in \mathcal{P}^{\text{co}}_{\tilde{\mathcal{E}}}}\sum_{\omega \in \Omega} f(\omega) {P}_{\tilde{\mathcal{E}}}(\{\omega\})&=\int_{\Omega} f \text{ d}\underline{P}_{\tilde{\mathcal{E}}}=\inf_{P_{\tilde{\mathcal{E}}} \in \mathcal{P}^{\text{co}}_{\tilde{\mathcal{E}}}}\int_{\Omega} f \text{ d}{P}_{\tilde{\mathcal{E}}}\\&=\inf_{P_{\tilde{\mathcal{E}}} \in \mathcal{P}^{\text{co}}_{\tilde{\mathcal{E}}}}\int_{\Omega} \hat{f} \text{ d}{P}_{\tilde{\mathcal{E}}}=\int_{\Omega} \hat{f} \text{ d}\underline{P}_{\tilde{\mathcal{E}}}=\inf_{P_{\tilde{\mathcal{E}}} \in \mathcal{P}^{\text{co}}_{\tilde{\mathcal{E}}}}\sum_{\omega \in \Omega} \hat{f}(\omega) {P}_{\tilde{\mathcal{E}}}(\{\omega\}),
\end{align*}
concluding the proof of (2). Notice also that 
\begin{align*}
\sup_{P_{\tilde{\mathcal{E}}} \in \mathcal{P}^{\text{co}}_{\tilde{\mathcal{E}}}}\sum_{\omega \in \Omega} f(\omega) {P}_{\tilde{\mathcal{E}}}(\{\omega\})&=\int_{\Omega} f \text{ d}\overline{P}_{\tilde{\mathcal{E}}}=\sup_{P_{\tilde{\mathcal{E}}} \in \mathcal{P}^{\text{co}}_{\tilde{\mathcal{E}}}}\int_{\Omega} f \text{ d}{P}_{\tilde{\mathcal{E}}}\\&=\sup_{P_{\tilde{\mathcal{E}}} \in \mathcal{P}^{\text{co}}_{\tilde{\mathcal{E}}}}\int_{\Omega} \hat{f} \text{ d}{P}_{\tilde{\mathcal{E}}}=\int_{\Omega} \hat{f} \text{ d}\overline{P}_{\tilde{\mathcal{E}}}=\sup_{P_{\tilde{\mathcal{E}}} \in \mathcal{P}^{\text{co}}_{\tilde{\mathcal{E}}}}\sum_{\omega \in \Omega} \hat{f}(\omega) {P}_{\tilde{\mathcal{E}}}(\{\omega\}).
\end{align*}

(3). By Theorem \ref{subadd_ergo}.(3) and the proof of claim (2), claim (3) follows.
\end{proof}

\begin{proof}[Proof of Theorem \ref{slln}]
%This proof is very similar to the one of \cite[Theorem 4]{cerreia}. The main difference is that we expressly request that $\Omega$ is finite or countable. 
By assumption, we have that $\mathbf{f}$ is stationary. Then, by a mathematical induction argument, we have that for all $k \in \mathbb{N}$ and all Borel subset $B$ of $\mathbb{R}$,
\begin{equation}\label{station}
\underline{P}_{\tilde{\mathcal{E}}}(\{\omega \in \Omega : f_1(\omega) \in B\})=\underline{P}_{\tilde{\mathcal{E}}}(\{\omega \in \Omega : f_k(\omega) \in B\}).
\end{equation} 
Equation \eqref{station} implies that for all $k \in \mathbb{N}$ and all Borel subset $B$ of $\mathbb{R}$,
$$\underline{P}_{\tilde{\mathcal{E}}}^{\mathbf{f}} (\{x \in \mathbb{R}^\mathbb{N}:x_k \in B\})=\underline{P}_{\tilde{\mathcal{E}}}(\{\omega \in \Omega : f_k(\omega) \in B\})=\underline{P}_{\tilde{\mathcal{E}}}(\{\omega \in \Omega : f_1(\omega) \in B\}).$$
Now, since $(f_k) \subset B(\Omega,\mathcal{F})$, we have that there exists $m \in \mathbb{R}$ such that $-m\mathbbm{1}_\Omega \leq f_1 \leq m\mathbbm{1}_\Omega$. By replacing $B$ with $[-m,m]$, we obtain that 
\begin{equation}\label{B9}
\underline{P}_{\tilde{\mathcal{E}}}^{\mathbf{f}} (\{x \in \mathbb{R}^\mathbb{N}:x_k \in [-m,m]\})=\underline{P}_{\tilde{\mathcal{E}}}(\{\omega \in \Omega : f_1(\omega) \in [-m,m]\})=1, \quad \text{for all } k \in \mathbb{N}.
\end{equation}
Let us now define the function $\pi:\mathbb{R}^\mathbb{N} \rightarrow \mathbb{R}$ as
$$x \mapsto \pi(x):=\begin{cases}
x_1 &\text{if } x_1 \in [-m,m]\\
0 &\text{otherwise}
\end{cases}.$$
We immediately see that $\pi$ belongs to $B(\mathbb{R}^\mathbb{N},\sigma(\mathcal{C}))$. Notice also that 
\begin{equation}\label{B10}
\bigcap_{k=1}^\infty \left\lbrace{x \in \mathbb{R}^\mathbb{N} : x_k \in [-m,m]}\right\rbrace \subset \bigcap_{k=1}^\infty \left\lbrace{x \in \mathbb{R}^\mathbb{N} : \frac{1}{k} \sum_{j=1}^k \pi(s^{j-1}(x))=\frac{1}{k} \sum_{j=1}^k x_j}\right\rbrace.
\end{equation}
Given \eqref{B9} and \eqref{B10}, and since $\underline{P}_{\tilde{\mathcal{E}}}^{\mathbf{f}}$ is both convex and continuous at $\mathbb{R}^\mathbb{N}$, we have that
\begin{equation}\label{B11}
\underline{P}_{\tilde{\mathcal{E}}}^{\mathbf{f}} \left( \bigcap_{k=1}^\infty \left\lbrace{x \in \mathbb{R}^\mathbb{N} : \frac{1}{k} \sum_{j=1}^k \pi(s^{j-1}(x))=\frac{1}{k} \sum_{j=1}^k x_j}\right\rbrace \right)=1.
\end{equation}
By \cite[Theorem 2]{cerreia} and the fact that $\underline{P}_{\tilde{\mathcal{E}}}^{\mathbf{f}}$ is shift invariant and ergodic, then there exists $\pi^* \in B(\mathbb{R}^\mathbb{N},\mathcal{G})$ such that 
\begin{equation}\label{B12}
\underline{P}_{\tilde{\mathcal{E}}}^{\mathbf{f}} \left( \left\lbrace{x \in \mathbb{R}^\mathbb{N} : \int_{\mathbb{R}^\mathbb{N}} \pi^\star \text{ d}\underline{P}_{\tilde{\mathcal{E}}}^{\mathbf{f}} \leq \lim_{k \rightarrow \infty} \frac{1}{k} \sum_{j=1}^k \pi(s^{j-1}(x))=\pi^\star(x) \leq \int_{\mathbb{R}^\mathbb{N}} \pi^\star \text{ d}\overline{P}_{\tilde{\mathcal{E}}}^{\mathbf{f}}}\right\rbrace \right)=1.
\end{equation}
Then, by \eqref{B11}, \eqref{B12}, and the fact that $\underline{P}_{\tilde{\mathcal{E}}}^{\mathbf{f}}$ is convex, we have that
\begin{equation}\label{B13}
\underline{P}_{\tilde{\mathcal{E}}}^{\mathbf{f}} \left( \left\lbrace{x \in \mathbb{R}^\mathbb{N} : \int_{\mathbb{R}^\mathbb{N}} \pi^\star \text{ d}\underline{P}_{\tilde{\mathcal{E}}}^{\mathbf{f}} \leq \lim_{k \rightarrow \infty} \frac{1}{k} \sum_{j=1}^k x_j=\pi^\star(x) \leq \int_{\mathbb{R}^\mathbb{N}} \pi^\star \text{ d}\overline{P}_{\tilde{\mathcal{E}}}^{\mathbf{f}}}\right\rbrace \right)=1.
\end{equation}
Define now the set $E:=\{x \in \mathbb{R}^\mathbb{N} : \lim_{k \rightarrow \infty} \frac{1}{k} \sum_{j=1}^k \pi(s^{j-1}(x))=\pi^\star(x)\}$ and the function $\pi_k:=\frac{1}{k} \sum_{j=1}^k \pi(s^{j-1})$, for all $k \in \mathbb{N}$. Then, by \eqref{B12} we have that $P(E)=1$, for all $P \in$ core$(\underline{P}_{\tilde{\mathcal{E}}}^{\mathbf{f}})$. By construction, $(\mathbbm{1}_E\pi_k)_{k \in \mathbb{N}} \subset B(\mathbb{R}^\mathbb{N},\sigma(\mathcal{C}))$ is uniformly bounded and converges (pointwise) to $\mathbbm{1}_E\pi^\star$. By $\underline{P}_{\tilde{\mathcal{E}}}^{\mathbf{f}}$ being convex, $P(E)$ being $1$ for all $P \in$ core$(\underline{P}_{\tilde{\mathcal{E}}}^{\mathbf{f}})$, and \cite[Theorem 22]{cerreia4}, it follows that that
\begin{align}\label{B14}
\begin{split}
\int_{\mathbb{R}^\mathbb{N}} \pi^\star \text{ d}\underline{P}_{\tilde{\mathcal{E}}}^{\mathbf{f}}&=\int_{\mathbb{R}^\mathbb{N}} \mathbbm{1}_E\pi^\star \text{ d}\underline{P}_{\tilde{\mathcal{E}}}^{\mathbf{f}}=\int_{\mathbb{R}^\mathbb{N}} \lim_{k \rightarrow \infty} \mathbbm{1}_E\pi_k \text{ d}\underline{P}_{\tilde{\mathcal{E}}}^{\mathbf{f}}\\&=\lim_{k \rightarrow \infty} \int_{\mathbb{R}^\mathbb{N}}  \mathbbm{1}_E\pi_k \text{ d}\underline{P}_{\tilde{\mathcal{E}}}^{\mathbf{f}}=\lim_{k \rightarrow \infty} \int_{\mathbb{R}^\mathbb{N}}  \pi_k \text{ d}\underline{P}_{\tilde{\mathcal{E}}}^{\mathbf{f}}.
\end{split}
\end{align}
Then, because $\underline{P}_{\tilde{\mathcal{E}}}^{\mathbf{f}}$ is convex and shift invariant, we have that, for all $n \in \mathbb{N}$,
\begin{align*}
\int_{\mathbb{R}^\mathbb{N}} \pi_k \text{ d}\underline{P}_{\tilde{\mathcal{E}}}^{\mathbf{f}} =\int_{\mathbb{R}^\mathbb{N}} \frac{1}{k} \sum_{j=1}^k \pi(s^{j-1}) \text{ d}\underline{P}_{\tilde{\mathcal{E}}}^{\mathbf{f}} \geq \frac{1}{k} \sum_{j=1}^k \int_{\mathbb{R}^\mathbb{N}} \pi(s^{j-1}) \text{ d}\underline{P}_{\tilde{\mathcal{E}}}^{\mathbf{f}} = \int_{\mathbb{R}^\mathbb{N}} \pi \text{ d}\underline{P}_{\tilde{\mathcal{E}}}^{\mathbf{f}}.
\end{align*}
Then, by \eqref{B14}, we have that $\int_{\mathbb{R}^\mathbb{N}} \pi^\star \text{ d}\underline{P}_{\tilde{\mathcal{E}}}^{\mathbf{f}} \geq \int_{\mathbb{R}^\mathbb{N}} \pi \text{ d}\underline{P}_{\tilde{\mathcal{E}}}^{\mathbf{f}}$. A similar argument shows that $\int_{\mathbb{R}^\mathbb{N}} \pi^\star \text{ d}\overline{P}_{\tilde{\mathcal{E}}}^{\mathbf{f}} \leq \int_{\mathbb{R}^\mathbb{N}} \pi \text{ d}\overline{P}_{\tilde{\mathcal{E}}}^{\mathbf{f}}$. Now, since by construction
$$\int_{\mathbb{R}^\mathbb{N}} \pi \text{ d}\underline{P}_{\tilde{\mathcal{E}}}^{\mathbf{f}} =\int_\Omega f_1 \text{ d} \underline{P}_{\tilde{\mathcal{E}}} \quad \text{and} \quad \int_{\mathbb{R}^\mathbb{N}} \pi \text{ d}\overline{P}_{\tilde{\mathcal{E}}}^{\mathbf{f}} =\int_\Omega f_1 \text{ d} \overline{P}_{\tilde{\mathcal{E}}},$$
we have that \eqref{B13} gives us
\begin{align*}
1&=\underline{P}_{\tilde{\mathcal{E}}}^{\mathbf{f}} \left( \left\lbrace{x \in \mathbb{R}^\mathbb{N} : \int_{\mathbb{R}^\mathbb{N}} \pi \text{ d}\underline{P}_{\tilde{\mathcal{E}}}^{\mathbf{f}} \leq \lim_{k \rightarrow \infty} \frac{1}{k} \sum_{j=1}^k x_j \leq \int_{\mathbb{R}^\mathbb{N}} \pi \text{ d}\overline{P}_{\tilde{\mathcal{E}}}^{\mathbf{f}}}\right\rbrace \right)\\
&=\underline{P}_{\tilde{\mathcal{E}}} \left( \left\lbrace{\omega \in \Omega: \int_\Omega f_1 \text{ d} \underline{P}_{\tilde{\mathcal{E}}} \leq \lim\limits_{k \rightarrow \infty} \frac{1}{k} \sum\limits_{j=1}^k f_j(\omega) \leq \int_\Omega f_1 \text{ d} \overline{P}_{\tilde{\mathcal{E}}}}\right\rbrace \right).
\end{align*}
Because $\underline{P}_{\tilde{\mathcal{E}}}$ is convex, by \cite[Theorem 38]{marinacci_ambig} we have that 
\begin{equation*}
\inf_{P_{\tilde{\mathcal{E}}} \in \mathcal{P}^{\text{co}}_{\tilde{\mathcal{E}}}}\sum_{\omega \in \Omega} f_1(\omega) {P}_{\tilde{\mathcal{E}}}(\{\omega\})= \inf_{P_{\tilde{\mathcal{E}}} \in \mathcal{P}^{\text{co}}_{\tilde{\mathcal{E}}}} \int_{\Omega} f_1 \text{ d}P =  \int_{\Omega} f_1 \text{ d}\underline{P}_{\tilde{\mathcal{E}}},
\end{equation*}
where the first equality comes from $\Omega$ being at most countable, and similarly
\begin{equation*}
\sup_{P_{\tilde{\mathcal{E}}} \in \mathcal{P}^{\text{co}}_{\tilde{\mathcal{E}}}}\sum_{\omega \in \Omega} f_1(\omega) {P}_{\tilde{\mathcal{E}}}(\{\omega\})= \sup_{P_{\tilde{\mathcal{E}}} \in \mathcal{P}^{\text{co}}_{\tilde{\mathcal{E}}}} \int_{\Omega} f_1 \text{ d}P =  \int_{\Omega} f_1 \text{ d}\overline{P}_{\tilde{\mathcal{E}}}.
\end{equation*}
Since $\underline{P}_{\tilde{\mathcal{E}}}$ is a lower probability, this implies that
$$\underline{P}_{\tilde{\mathcal{E}}} \left( \left\lbrace{\omega \in \Omega: \inf_{P_{\tilde{\mathcal{E}}} \in \mathcal{P}^{\text{co}}_{\tilde{\mathcal{E}}}}\sum_{\omega \in \Omega} f_1(\omega) {P}_{\tilde{\mathcal{E}}}(\{\omega\}) \leq \lim\limits_{k \rightarrow \infty} \frac{1}{k} \sum\limits_{j=1}^k f_j(\omega) \leq \sup_{P_{\tilde{\mathcal{E}}} \in \mathcal{P}^{\text{co}}_{\tilde{\mathcal{E}}}}\sum_{\omega \in \Omega} f_1(\omega) {P}_{\tilde{\mathcal{E}}}(\{\omega\})}\right\rbrace \right)$$
is equal to $1$, which proves the statement.
\end{proof}

\begin{proof}[Proof of Lemma \ref{suff_cond_invariance}]
Given our assumptions, we have that for any $A^\prime\in\mathcal{F}$ we can always find $P^{A^\prime}_{\mathcal{E}_t} \in \{P^A_{\mathcal{E}_t}\}_{A\in\mathcal{F}}$ such that 
$$\underline{P}_{\mathcal{E}_t}(A^\prime)=P^{A^\prime}_{\mathcal{E}_t}(A^\prime)=P^{A^\prime}_{\mathcal{E}_t}(T^{-1}(A^\prime))=\underline{P}_{\mathcal{E}_t}(T^{-1}(A^\prime)),$$
so $\underline{P}_{\mathcal{E}_t}$ is $T$-invariant. Because this holds for all $t\geq T$, it also holds for a collection $\{P^A_{\tilde{\mathcal{E}}}\}_{A\in\mathcal{F}}$ belonging to the almost sure limit $\mathcal{P}^{\text{co}}_{\tilde{\mathcal{E}}}$ of sequence $(\mathcal{P}^{\text{co}}_{\mathcal{E}_t})$. In turn, this implies that $\underline{P}_{\tilde{\mathcal{E}}}$ is $T$-invariant.
\end{proof}

\begin{proof}[Proof of Lemma \ref{suff_cond_inv_ergo}]
Suppose that there exist $T\in\mathbb{N}_0$ and  $P_{\mathcal{E}_T}^\prime\in\mathcal{P}^{\text{co}}_{\mathcal{E}_T}$ such that $P_{\mathcal{E}_T}^\prime(A)=0$, for all $A\in\mathcal{G}$. This implies that $\underline{P}_{\mathcal{E}_T}(A)=0$. Then, we have that 
\begin{align*}
    P_{\mathcal{E}_{T+1}}^\prime(A)&=\sum_{\emptyset\neq E\in\mathcal{E}_{T+1}} \frac{P_{\mathcal{E}_T}^\prime(A\cap E)}{P_{\mathcal{E}_T}^\prime(E)}P_{\mathcal{E}_{T+1}}^\prime(E)\\
    &\leq \sum_{\emptyset\neq E\in\mathcal{E}_{T+1}} \frac{P_{\mathcal{E}_T}^\prime(A)}{P_{\mathcal{E}_T}^\prime(E)}P_{\mathcal{E}_{T+1}}^\prime(E)=0.
\end{align*}
So $P_{\mathcal{E}_{T+1}}^\prime(A)=0$, which implies $\underline{P}_{\mathcal{E}_{T+1}}(A)=0$. A similar argument shows that $P_{\mathcal{E}_{t}}^\prime(A)=0$, for all $t \geq T$, which implies that $P_{\tilde{\mathcal{E}}}^\prime(A)=0$. In turn, this implies that $\underline{P}_{\tilde{\mathcal{E}}}(A)=0$. But because this is true for all $A\in\mathcal{G}$, we have that $\underline{P}_{\tilde{\mathcal{E}}}$ is ergodic.
\end{proof}

\begin{proof}[Proof of Lemma \ref{suff_cond_stricter_than_before}]
Pick any $A,B\in\mathcal{F}$ such that $A\subset B$. If (i) holds, it is immediate to see that for all $t\geq T$, $P_{\mathcal{E}_{t+1}}^\prime(A)\leq P_{\mathcal{E}_{t+1}}(A)$ and $P_{\mathcal{E}_{t+1}}^\prime(B)\leq P_{\mathcal{E}_{t+1}}(B)$, for all $P_{\mathcal{E}_{t+1}}\in \mathcal{P}^{\text{co}}_{\mathcal{E}_{t+1}}$. This implies that $P_{\mathcal{E}_{t+1}}^\prime(A)= \underline{P}_{\mathcal{E}_{t+1}}(A)$ and $P_{\mathcal{E}_{t+1}}^\prime(B)= \underline{P}_{\mathcal{E}_{t+1}}(B)$. In turn, because this holds for all $t\geq T$, we have that the almost sure limit $P_{\tilde{\mathcal{E}}}^\prime$ of $P_{\mathcal{E}_{t+1}}^\prime$ is such that $P_{\tilde{\mathcal{E}}}^\prime(A)= \underline{P}_{\tilde{\mathcal{E}}}(A)$ and $P_{\tilde{\mathcal{E}}}^\prime(B)= \underline{P}_{\tilde{\mathcal{E}}}(B)$. This implies that $\underline{P}_{\tilde{\mathcal{E}}}$ is convex by \cite[Theorem 38.(ii)]{marinacci_ambig} and \cite{walley_report}.

Pick any finite chain $(A_i)_{i=1}^n \subset \mathcal{F}$. If (ii) holds, it is immediate to see that for all $t\geq T$, $P_{\mathcal{E}_{t+1}}^\prime(A_i)\leq P_{\mathcal{E}_{t+1}}(A_i)$, for all $i\in\{1,\ldots,n\}$. This implies that $P_{\mathcal{E}_{t+1}}^\prime(A_i)= \underline{P}_{\mathcal{E}_{t+1}}(A_i)$, for all $i\in\{1,\ldots,n\}$. In turn, because this holds for all $t\geq T$, we have that the almost sure limit $P_{\tilde{\mathcal{E}}}^\prime$ of $P_{\mathcal{E}_{t+1}}^\prime$ is such that $P_{\tilde{\mathcal{E}}}^\prime(A_i)= \underline{P}_{\tilde{\mathcal{E}}}(A_i)$, for all $i\in\{1,\ldots,n\}$. This implies that $\underline{P}_{\tilde{\mathcal{E}}}$ is convex by \cite[Theorem 38.(iii)]{marinacci_ambig} and \cite{walley_report}. A similar argument combined with \cite[Theorem 38.(iv)]{marinacci_ambig} and \cite[Corollary 2]{marinacci_napol} gives us that condition (iii) implies $\underline{P}_{\tilde{\mathcal{E}}}$ being convex.
\end{proof}

\begin{proof}[Proof of Lemma \ref{suff_cond_stricter_than_before2}]
We prove the three points separately.
\begin{itemize}
    \item[(i)] By \cite[Theorem 10]{marinacci_ambig}, we have that if
    \begin{itemize}
        \item[\textbullet] core$(\underline{P}_{\tilde{\mathcal{E}}})=\mathcal{P}^{\text{co}}_{\tilde{\mathcal{E}}}$ is nonempty,
        \item[\textbullet] $\underline{P}_{\tilde{\mathcal{E}}}(A)=\min_{P\in\mathcal{P}^{\text{co}}_{\tilde{\mathcal{E}}}}P(A)$, for all $A\in\mathcal{F}$,
        \item[\textbullet] $\mathcal{P}^{\text{co}}_{\tilde{\mathcal{E}}}$ is a weakly compact subset of the space $ca(\mathcal{F})$ of all measures on $\mathcal{F}$ having finite total variation norm,
    \end{itemize}
       then $\underline{P}_{\tilde{\mathcal{E}}}$ is continuous at $\Omega$. The first condition is always satisfied by $\underline{P}_{\tilde{\mathcal{E}}}$ by construction. The second condition is always satisfied by $\underline{P}_{\tilde{\mathcal{E}}}$ by construction and by the fact that core$(\underline{P}_{\tilde{\mathcal{E}}})=\mathcal{P}^{\text{co}}_{\tilde{\mathcal{E}}}$ is weak$^\star$-compact. The last condition is always satisfied thanks to \cite[Lemma 9]{marinacci_ambig}. In turn, we conclude that $\underline{P}_{\tilde{\mathcal{E}}}$ is always continuous at $\Omega$.
    \item[(ii)] Immediate from our assumption and the fact that $(\mathcal{P}^{\text{co}}_{\mathcal{E}_t})$ converges (almost surely) to $\mathcal{P}^{\text{co}}_{\tilde{\mathcal{E}}}$.
    \item[(iii)] Recall that the DPK update of probability measure $P$ given partition $\mathcal{E}$ of $\Omega$ is given by $P_\mathcal{E}=\sum_{\emptyset\neq E\in\mathcal{E}}P(A\mid E) [\beta(n) P(E)+(1-\beta(n)) P^{emp}(E)]$, where $\beta(n)$ is a coefficient in $[0,1]$ depending on the amount $n$ of data available and chosen subjectively by the agent, and $P^{emp}$ is a well defined probability measure \cite[Section 4]{prob.kin}.\footnote{As the updating process continues, the amount of data available and probability function $P^{emp}$ need to be indexed by time $t$, so we write $n_t$ and $P^{emp}_t$ \cite[Section 5]{prob.kin}.} If the assumptions in (iii) hold, we have that
    \begin{align*}
        P_{\mathcal{E}_{t+1}}(A)&=\sum_{\emptyset\neq E\in\mathcal{E}_{t+1}} \frac{P_{\mathcal{E}_t}(A\cap E)}{P_{\mathcal{E}_t}(E)} \left[\beta(n_{t+1}) P_{\mathcal{E}_t}(E)+(1-\beta(n_{t+1})) P^{emp}_{t+1}(E) \right]\\
        &=\sum_{\emptyset\neq E\in\mathcal{E}_{t+1}} \frac{P_{\mathcal{E}_t}(T^{-1}(A)\cap E)}{P_{\mathcal{E}_t}(E)} \left[\beta(n_{t+1}) P_{\mathcal{E}_t}(E)+(1-\beta(n_{t+1})) P^{emp}_{t+1}(E) \right]\\
        &=P_{\mathcal{E}_{t+1}}(T^{-1}(A)),
    \end{align*}
    for all $A\in\mathcal{F}$ and all $P_{\mathcal{E}_t} \in \mathcal{P}^{\text{co}}_{\mathcal{E}_t}$. This implies that $\mathcal{P}^{\text{co}}_{\mathcal{E}_{t+1}} \subset \mathcal{I}$. Because this is true for all $t\geq T$, we have that the almost sure limit $\mathcal{P}^{\text{co}}_{\tilde{\mathcal{E}}}$ too is a subset of $\mathcal{I}$, which implies that $\underline{P}_{\tilde{\mathcal{E}}}$ is functionally invariant.
\end{itemize}
\end{proof}

\begin{proof}[Proof of Lemma \ref{lemma1}]
We know that $(S_k)$ satisfies \eqref{bd_seq}, so $(S_k) \subset B(\Omega,\mathcal{F})$. If $(S_k)$ is superadditive and $\mathcal{P}^{\text{co}}_{\tilde{\mathcal{E}}} \subset \mathcal{I}$, then, for all $k,\ell \in \mathbb{N}$, we have the following
\begin{align*}
-a_{k+\ell}&=\inf_{P_{\tilde{\mathcal{E}}}\in\mathcal{P}^{\text{co}}_{\tilde{\mathcal{E}}}} \sum_{\omega \in \Omega} S_{k +\ell}(\omega)P_{\tilde{\mathcal{E}}}(\{\omega\}) \geq  \inf_{P_{\tilde{\mathcal{E}}}\in\mathcal{P}^{\text{co}}_{\tilde{\mathcal{E}}}}\sum_{\omega \in \Omega} \left[ S_{k}(\omega) + (S_\ell \circ T^k)(\omega) \right]P_{\tilde{\mathcal{E}}}(\{\omega\})\\
&\geq  \inf_{P_{\tilde{\mathcal{E}}}\in\mathcal{P}^{\text{co}}_{\tilde{\mathcal{E}}}}\sum_{\omega \in \Omega} S_{k}(\omega) P_{\tilde{\mathcal{E}}}(\{\omega\}) +  \inf_{P_{\tilde{\mathcal{E}}}\in\mathcal{P}^{\text{co}}_{\tilde{\mathcal{E}}}}\sum_{\omega \in \Omega} (S_\ell \circ T^k)(\omega) P_{\tilde{\mathcal{E}}}(\{\omega\})\\
&= \inf_{P_{\tilde{\mathcal{E}}}\in\mathcal{P}^{\text{co}}_{\tilde{\mathcal{E}}}}\sum_{\omega \in \Omega} S_{k}(\omega) P_{\tilde{\mathcal{E}}}(\{\omega\}) +  \inf_{P_{\tilde{\mathcal{E}}}\in\mathcal{P}^{\text{co}}_{\tilde{\mathcal{E}}}}\sum_{\omega \in \Omega} S_\ell (\omega) P_{\tilde{\mathcal{E}}}(\{\omega\})=-a_k-a_\ell.
\end{align*}
A similar procedure shows the result when $(S_k)$ is subadditive and $a_k$ is the supremum over $\mathcal{P}^{\text{co}}_{\tilde{\mathcal{E}}}$ of the sum of $S_{k}(\omega)P_{\tilde{\mathcal{E}}}(\{\omega\})$.
\end{proof}

\begin{proof}[Proof of Lemma \ref{core_conv}]
To prove the claim, we show that if we consider a countably additive probability $P^\prime \in \Delta(\Omega,\mathcal{F})$ that \textit{does not} belong to $\text{Conv}(P_1,\ldots,P_k)$, then there exists $X^\prime \in B(\Omega,\mathcal{F})$ such that either $P^\prime(X^\prime) < \underline{P}(X)$ or $P^\prime(X^\prime) > \overline{P}(X)$.

For every element $X\in B(\Omega,\mathcal{F})$, let $P_L(X):=\inf_{P\in\text{Conv}(P_1,\ldots,P_k)} P(X)$ and also  $P_U(X):=\sup_{P\in\text{Conv}(P_1,\ldots,P_k)} P(X)$. So $[P_L(X),P_U(X)]$ is the closed interval of previsions for $X$ with respect to the convex hull $\text{Conv}(P_1,\ldots,P_k)$. Then, by the Krein-Milman theorem, we know that $\text{Conv}(P_1,\ldots,P_k)$ is the largest (closed) convex set satisfying these $\#B(\Omega,\mathcal{F})$ many pairs of constraints. Since $P^\prime \not\in \text{Conv}(P_1,\ldots,P_k)$, it fails at least one of these $\#B(\Omega,\mathcal{F})$ many pairs of constraints, so $P^\prime \not\in \text{core}_B(\underline{P})$.
\end{proof}

\begin{proof}[Proof of Lemma \ref{lemma_marinacci_nice}]
It is immediate to see that $\underline{P}_{\tilde{\mathcal{E}}}$ is a positive exact game on $\mathcal{F}$. In addition, we know by Lemma \ref{suff_cond_stricter_than_before2}.(i) that  $\underline{P}_{\tilde{\mathcal{E}}}$ is always continuous at $\Omega$. This implies by Claim \ref{claim_imp_marin} that $\underline{P}_{\tilde{\mathcal{E}}}$ is continuous at each $A$, that is,  $\underline{P}_{\tilde{\mathcal{E}}}$ is continuous in the sense of section \ref{NAD_imp}.(iii).
\end{proof}

\bibliographystyle{plain}
\bibliography{ergodic_theory} 
\end{document}